\documentclass[11pt]{article}
\usepackage{amsmath}
\usepackage{amssymb}
\usepackage{amscd}
\usepackage{latexsym}
\usepackage{amsthm}
\usepackage{epsfig}
\usepackage{setspace}
\usepackage{amsfonts}
\usepackage[all]{xypic}
 \usepackage{epstopdf,ifpdf}
\ifpdf
  \usepackage{pdfsync}
  \usepackage{hyperref}
\fi

\newtheorem{theorem}{Theorem}[section]
\newtheorem{lemma}[theorem]{Lemma}
\newtheorem{proposition}[theorem]{Proposition}
\newtheorem{corollary}[theorem]{Corollary}

{
\theoremstyle{definition}
\newtheorem{definition}[theorem]{Definition}
\newtheorem{example}[theorem]{Example}

\newtheorem{remark}[theorem]{Remark}

}

\newcommand{\nc}{\newcommand}
\newcommand{\Spec}{\operatorname{Spec}}

\newcommand{\Span}{\operatorname{Span}}
\newcommand{\rank}{\operatorname{rank}}

\renewcommand{\Im}{\operatorname{Im}}

\newcommand{\id}{\operatorname{id}}

\renewcommand{\dim}{\operatorname{dim}}

\newcommand{\Sym}{\operatorname{Sym}}

\newcommand{\gr}{\operatorname{gr}}

\newcommand{\rad}{\operatorname{rad}}
\newcommand{\gmod}{\mathsf{gmod}}

\newcommand{\Z}{\mathbb{Z}}
\newcommand{\Q}{\mathbb{Q}}

\newcommand{\R}{\mathbb{R}}
\newcommand{\C}{\mathbb{C}}
\newcommand{\V}{\mathbb{V}}

\newcommand{\lra}{\longrightarrow}

\newcommand{\Lplusmu}{\Lambda^+_\mu}

\newcommand{\la}{\lambda}

\newcommand{\subs}{\subseteq}

\renewcommand{\a}{\alpha}
\renewcommand{\b}{\beta}
\newcommand{\surj}{\twoheadrightarrow}
\newcommand{\Hom}{\operatorname{Hom}}
\newcommand{\End}{\operatorname{End}}

\newcommand{\cI}{\mathcal{I}}
\newcommand{\cZ}{\mathcal{Z}}
\newcommand{\cJ}{\mathcal{J}}
\newcommand{\z}{\zeta}
\newcommand{\ga}{\gamma}
\nc{\ep}{\varepsilon}
\nc{\excise}[1]{}
\renewcommand{\t}{\mathfrak{t}}
\newcommand{\MV}{\mathfrak{M}(\cV)}
\newcommand{\cV}{\mathcal{V}}
\renewcommand{\cH}{\mathcal{H}}
\newcommand{\mh}{\mathfrak{h}}
\newcommand{\mt}{\mathfrak{t}}
\newcommand{\mb}{\mathfrak{b}}
\newcommand{\mg}{\mathfrak{g}}
\renewcommand{\mp}{\mathfrak{p}}
\newcommand{\mm}{\mathfrak{m}}

\newcommand{\cO}{\mathcal{O}}

\newcommand{\fg}{\mathfrak{g}}

\newcommand{\wt}{\widetilde}
\newcommand{\wh}{\widehat}

\newcommand{\Lie}{\operatorname{Lie}}

\newcommand{\IH}{\operatorname{IH}}

\newcommand{\Fl}{X}

\begin{document}

\noindent {\Large \bf Localization algebras and deformations of Koszul algebras}
\bigskip\\
{\bf Tom Braden}\footnote{Supported by NSA grant H98230-08-1-0097.}\\
Department of Mathematics and Statistics, University of Massachusetts,
Amherst, MA 01003\smallskip \\
{\bf Anthony Licata}\\
Department of Mathematics, Stanford University,
Palo Alto, CA 94305\smallskip \\
{\bf Christopher Phan}\\
Department of Mathematics, Bucknell University, 
Lewisburg, PA\smallskip \\
{\bf Nicholas Proudfoot}\footnote{Supported by NSF grants DMS-0738335 and DMS-0950383.}\\
Department of Mathematics, University of Oregon,
Eugene, OR 97403\smallskip \\
{\bf Ben Webster}\footnote{Supported by an NSF Postdoctoral Research Fellowship and by NSA grant H98230-10-1-0199.}\\
Department of Mathematics, University of Oregon,
Eugene, OR 97403
\bigskip\\
{\small
\begin{quote}
\noindent {\em Abstract.}
We show that the center of a flat graded deformation of a standard Koszul algebra
$A$ behaves in many ways like the torus-equivariant cohomology ring of 
an algebraic variety with finite fixed-point set.  
In particular, the center of $A$ acts by characters on the deformed standard modules,
providing a ``localization map."  We construct a universal graded
deformation of $A$, and show that the spectrum of its center is 
supported on a certain arrangement of hyperplanes which
is orthogonal to the arrangement coming from the algebra Koszul dual to $A$.
This is an algebraic version of a duality discovered by
Goresky and MacPherson between the equivariant cohomology
rings of partial flag varieties and Springer fibers; we 
recover and generalize their result by showing that the
center of the universal deformation for the ring governing 
a block of parabolic category $\cO$ for $\mathfrak{gl}_n$ is 
isomorphic to the equivariant cohomology of a Spaltenstein 
variety.  We also identify the center of the deformed 
version of the ``category $\cO$" of a hyperplane arrangement
(defined by the authors in a previous paper) with the equivariant cohomology of
a hypertoric variety.
\end{quote}
}
\bigskip
\begin{spacing}{1.2}

\section{Introduction}
\label{sec:introduction}
In 1976, Bernstein, Gelfand, and Gelfand introduced the category
$\cO$ of representations of a semisimple Lie algebra 
$\fg = \operatorname{Lie}(G)$ \cite{BGG}.  Over the course of the
next decade, several new techniques appeared in the algebraic and geometric study of this category.   Two of the most important were 
\begin{itemize}
\setlength{\itemsep}{-1pt}
\item the use of a deformed category $\widehat \cO$,
  which consists of families of representations over a formal neighborhood of $0$ in the weight space $\mh^*$ of $\fg$
\item connections to the geometry of the flag variety $G/B$,
  especially through the localization theorem of Beilinson and
  Bernstein \cite{BB}.
\end{itemize}
The first of these was used by Soergel \cite{Soe90} to
show that an integral block of category $\cO$ is equivalent to
the module category of a certain finite dimensional algebra $A$.
Furthermore, Soergel showed that the center of $A$ is isomorphic to the cohomology ring of $H^*(G/P)$, where $P\subset G$ is a parabolic
subalgebra that depends on the block.  This fact reflects a second connection
to geometry: each block of category $\cO$ not only has a geometric interpretation
via the localization theorem, it is also Koszul dual to the category of
Schubert smooth perverse sheaves on $G/P$, which is proved independently of the localization theorem in \cite{BGS96}.  
In later work \cite{Soe92}, Soergel showed that
$\widehat \cO$ can be described using the $T$-equivariant geometry of $G/P$,
where $T$ is a maximal torus of $G$.  In particular, he computed a deformation $\widehat A$ of $A$ whose module category is isomorphic to the corresponding block 
of $\widehat \cO$, and he showed that the center of $\widehat A$
is isomorphic to a completion of the equivariant cohomology ring $H^*_T(G/P)$.

Our aim in this paper is to study how the first of these techniques, that of studying the representation theory of a finite dimensional Koszul algebra by deforming it, can be applied in a general algebraic context, without the benefit of the geometric or Lie theoretic interpretations of category $\cO$. 
In Section \ref{sec:canon-deform-quadr} we
use the work of Braverman and Gaitsgory \cite{BG96} to show 
that any Koszul algebra $A$ has a universal flat graded deformation 
$\tilde{A}$, so that any other graded flat deformation is obtained from 
$\tilde{A}$ by a unique base change.
In Sections \ref{category O examples} and  \ref{sec:deformed O}, we use techniques of Soergel and Fiebig \cite{Soe90, Soe92, Fie03, Fie06, Fie08} 
to show that Soergel's deformed algebra is in fact the completion of the universal deformation.

For what follows, we will need to assume not just that our algebras are Koszul,
but also {\bf standard Koszul} (Definition \ref{def:stdk}), 
which roughly means that its module category is highest weight in the sense of \cite{CPS}.  The general study of standard Koszul algebras was initiated in \cite{ADL03}, and the main examples are the algebras $A$ introduced above, the generalizations obtained by replacing $\cO$ with its parabolic version, and
the algebras which were defined in \cite{GDKD} 
using hyperplane arrangements.

We focus our attention primarily on the center of the universal deformation $\tilde A$ 
of a standard Koszul algebra $A$.
Recall that Soergel identified the center of his deformed algebra
with the (completed) 
equivariant cohomology ring of a partial flag variety.
We generalize this result by showing that the center of $\tilde A$ always
``behaves like'' a torus-equivariant cohomology ring.
More precisely, we introduce a structure called a {\bf localization algebra}
(Definition \ref{gmalg}),
which is an abstraction of the data given by a torus-equivariant cohomology ring and the
localization map to the fixed point set, and we prove the following theorem
(Corollary \ref{gma}).

\begin{theorem}\label{first}
The center $\cZ(A) := Z(\tilde A)$ of the universal deformation
canonically admits the structure of a localization algebra.
\end{theorem}

For many of the examples of standard Koszul algebras mentioned above,
the localization algebras that we obtain are in fact isomorphic to
equivariant cohomology rings.
We have already addressed the case where $A$ is the algebra whose
module category is a block of $\cO$.
If $\fg = \mathfrak{gl}_n$ and $\cO$ is replaced by its parabolic version, 
then Brundan and Stroppel \cite{Bru06,Str09} show
that the center of $A$ is isomorphic to the cohomology ring of a Spaltenstein
variety.\footnote{
A Spaltenstein variety is a certain subvariety of $G/P$, where once again
the choice of $P$ depends on the choice of central character.  If the central character is generic,
then $P$ will be a Borel subalgebra, and the Spaltenstein variety will be a Springer fiber;
this is the case proven in \cite{Str09}.
See Section \ref{category O examples} for more details.}  In this case, we show
that the center of $\tilde A$ is isomorphic to the torus-equivariant cohomology
ring of the same variety (Theorem \ref{cat O and Spaltensteins}), generalizing
Soergel's result in the non-parabolic case and Brundan's result for the
un-deformed algebra.  Finally, in the case where $A$ is one of the
algebras introduced in \cite{GDKD}, we show that the center of $\tilde A$
is isomorphic to the cohomology ring of a hypertoric variety (Theorem \ref{eqb}),
generalizing the un-deformed, non-equivariant result of \cite[4.16]{GDKD}.

Once we have established that the center of the universal deformation
is a localization algebra (Theorem \ref{first}), we study 
the relationship between the localization algebras
associated to a dual pair of standard Koszul algebras.
This problem is motivated by an observation made by Goresky and MacPherson
in a paper that has, {\it a priori}, nothing to do with Koszul algebras \cite{GM}.\footnote{In
particular, it is unrelated to the appearance of Koszul duality in \cite{GKM},
which is very different in flavor from anything in this paper.} 
If $X$ is a variety equipped
with the action of a torus $T$ with isolated fixed points, then the localization map
in equivariant cohomology may be used to define a finite collection of vector subspaces
of $H^T_2(X)$, each of which is isomorphic to the Lie algebra $\mathfrak{t}$ of $T$.
Goresky and MacPherson observed that the arrangement associated to a partial flag variety for
$\mathfrak{gl}_n$ is dual to the one
associated to a Springer fiber, in the sense that the equivariant second homology groups
are dual as vector spaces, and the subspaces appearing on one side are the perpendicular
spaces to the subspaces appearing on the other side.

Our approach to the problem is to interpret and generalize the examples of \cite{GM} 
in a purely algebraic context.  First, we define what it means for two localization algebras
to be {\bf dual}, so that the result of \cite{GM} may be formulated as the duality of a certain pair of localization algebras.
We next introduce one more technical hypothesis:
we call an algebra $A$ {\bf flexible} if it is standard Koszul and the natural map
from the center of $\tilde A$ to the center of $A$ is surjective in degree 2.
Our main result (Theorem \ref{bang-flexible} 
and Corollary \ref{main}) is the following.

\begin{theorem}\label{second}
If $A$ is flexible, then so is the dual algebra $A^!$, and 
the localization algebras $\cZ(A)$ and $\cZ(A^!)$ are canonically dual.
\end{theorem}

\noindent
In light of Theorem \ref{cat O and Spaltensteins} and the fact 
that a regular block of parabolic category $\cO$ for $\mathfrak{gl}_n$
is Koszul dual to a singular block of ordinary category $\cO$ \cite{BGS96}, the examples
found by Goresky and MacPherson follow from Theorem \ref{second}.
In fact, a theorem of Backelin's \cite{Back99} says that an integral block of parabolic category $\cO$
is dual to another such block, 
and so the Goresky-MacPherson phenomenon generalizes to 
all $\mathfrak{gl}_n$ Spaltenstein varieties.

\begin{example}
Consider the quiver
$$
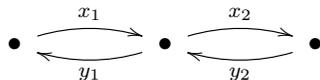
\xy
(-20,0)*++{\bullet}="1"; (0,0)*++{\bullet}="2";(20,0)*++{\bullet}="3";
{\ar@/^/^{x_1} "1";"2"}; {\ar@/^/^{y_1} "2";"1"};
{\ar@/^/^{x_2} "2";"3"}; {\ar@/^/^{y_2} "3";"2"};
\endxy$$
and let $A$ be the path algebra modulo the relations $x_1y_1$ and $y_1x_1 - x_2y_2$.
This is a noncommutative graded algebra, and it is standard Koszul.
The category of finitely generated $A$-modules 
is equivalent to a singular integral block of category $\cO$
for $\mathfrak{sl}_3$.  Then Soergel \cite{Soe90}
tells us that its center $Z(A)$ is isomorphic to the cohomology ring of $\mathbb{P}^2$.
Indeed, $Z(A)$ is generated by the degree 2 class $x_2y_2 + y_2x_2$, 
whose cube is zero.  The center $\cZ(A) := Z(\tilde A)$ of the universal deformation of $A$
is isomorphic to the $T^2$-equivariant cohomology ring of $\mathbb{P}^2$.

Now consider the Koszul dual ring $A^!$ of $A$, which is a quotient of the path algebra of
the dual quiver
$$\xy
(-20,0)*++{\bullet}="1"; (0,0)*++{\bullet}="2";(20,0)*++{\bullet}="3";
{\ar@/_/_{x_1^!} "2";"1"}; {\ar@/_/_{y_1^!} "1";"2"};
{\ar@/_/_{x_2^!} "3";"2"}; {\ar@/_/_{y_2^!} "2";"3"};
\endxy$$
by the relations $x_1^!y_1^! + y_2^!x_2^!$, $x_2^!y_2^!$, $x_2^!x_1^!$, and $y_1^!y_2^!$.
The module category of $A^!$ is equivalent to a regular block of parabolic category $\cO$
for the parabolic subalgebra of $\mathfrak{sl}_3$ preserving a line in $\C^3$ \cite[5.2.1]{Str03}.
The results of Brundan or Stroppel \cite{Bru06,Str09} tell us that the center of $A$ should be isomorphic
to the cohomology ring of a certain Springer fiber, namely the one consisting of two projective
lines that touch in a single point.
Indeed, the center of $A^!$ is generated by the degree 2 classes $y_1^!x_1^!$ and $y_2^!x_2^!$, 
with all products trivial. 
The center $\cZ(A^!) := Z(\tilde A^!)$ of the universal deformation of $A^!$
is isomorphic to the $T^1$-equivariant cohomology ring of the same variety.

The equivariant cohomology rings associated to these two algebras constitute the
simplest nontrivial example of a dual pair from \cite{GM}.  More details from the perspective
of \cite{GM} are given in Examples \ref{ptwo}, \ref{two lines}, and \ref{first duality}.
\end{example}

\begin{remark}
The following is meant only
to provide some additional geometric motivation for our results.
In both of the families of standard Koszul algebras considered in Sections \ref{hypertoric examples}-\ref{sec:deformed O}, the localization algebras that arise are isomorphic to equivariant
cohomology rings of certain algebraic symplectic manifolds or orbifolds.\footnote{Hypertoric
varieties are symplectic orbifolds, and Spaltenstein varieties are torus-equivariant
deformation retracts of resolved Slodowy slices, which are symplectic manifolds.}
We expect that the algebra itself will be isomorphic to  
the Ext-algebra of a certain module over a quantization of the structure sheaf of the manifold.
The map from the cohomology ring to the center of our algebra will then be induced
by the action of the constant sheaf on this module.
For the case studied in Section \ref{hypertoric examples}, this program is being carried
out in \cite{BLPWtorico}.  In the case of an algebra whose module category is equivalent to
an integral block of ordinary category $\cO$, it is well-known:
the manifold in question is $T^*(G/P)$,
and the module is the microlocalization of the direct sum of all of the Schubert-smooth
simple D-modules on $G/P$.  The authors plan to treat the case of parabolic
category $\cO$ in a forthcoming paper.

When two algebraic symplectic manifolds give rise to 
dual standard Koszul rings in this way, we refer to them as 
a {\bf symplectic dual} pair.  So the main result of 
this paper could be interpreted as saying that 
symplectic dual pairs have equivariant cohomology rings that are dual as localization algebras.
Beside the pairs of hypertoric varieties and pairs of resolved Slodowy slices
that we consider in this paper, other conjectural examples
include Hilbert schemes on ALE spaces, which we 
expect to be dual to certain moduli spaces of instantons on $\C^2$, and quiver varieties of simply
laced Dynkin type, which we expect to be dual to resolutions of slices to certain subvarieties
of the affine Grassmannian.
We expect further examples to arise from physics as Higgs branches of the moduli
space of vacua for mirror dual 3-dimensional $\mathcal{N}=4$ superconformal field theories,
or as the Higgs and Coulomb branches of a single such theory.  That 
hypertoric varieties occur in mirror dual theories was observed by Kapustin and Strassler in
\cite{KS99}.
\end{remark}

\subsection*{Acknowledgments}
The authors thank Jon Brundan, Andrew Connor, Peter Fiebig, Mark Goresky, 
Brad Shelton, and Catharina Stroppel for useful conversations, and the referee for 
helpful comments.  
T.B.\ thanks Reed College for its hospitality during the writing of this paper. 
A.L.\ thanks the Max Planck Institute for Mathematics in Bonn for hospitality.

\section{Localization algebras}\label{localg}

In this section we introduce localization algebras and define a notion of duality between them.  Before giving the general definition of a localization algebra, we consider a motivating setup from equivariant topology.

Let $X$ be a complex algebraic variety and $T$ be an algebraic torus acting on $X$ with the fixed point set $X^T$ finite and non-empty.  
The equivariant cohomology ring\footnote{All cohomology groups in this paper will be
taken with complex coefficients.} $H_T^*(X)$ is a graded algebra
over the polynomial ring $\Sym \t^*$ (the $T$-equivariant cohomology of a point), 
where the elements of $\t^*$ lie in degree 2.
From these data, we obtain natural graded algebra homomorphisms
\begin{equation}\label{algebras}
\Sym \t^* \hookrightarrow \Sym H^2_T(X) \to H^*_T(X)\to H^*_T(X^T) \cong H^*(X^T)\otimes \Sym\t^*,
\end{equation}
where the second map is given by multiplying classes of degree 2 together, and the third,
often called the {\bf localization map}, is
given by restriction to $X^T$.  
(To simplify notation, we will write $\otimes$ in place of $\otimes_\C$
throughout the paper.)
If the ordinary cohomology of $X$ vanishes in odd degrees (for example if $X$ is a smooth projective variety), then $H^*_T(X)$ is a free $\Sym\t^*$-module, the localization map is injective,
and the cokernel of the localization map is a torsion $\Sym\t^*$-module.
With this example in mind, we formulate the following general definitions.

\begin{definition}\label{gmalg}
A {\bf localization algebra}
is a quadruple $\cZ = (U,Z,\cI,h)$, where $U$ is a finite-dimensional
complex vector space, $Z$ is a finitely generated graded $\Sym U$-algebra, $\cI$ is a finite
set, and 
 $$h:Z\to \bigoplus_{\a\in\cI} \Sym U$$
is a homomorphism of $\Sym U$-algebras.
If the kernel and cokernel of $h$ are torsion $\Sym U$-modules, then we call $\cZ$ {\bf strong}.
If $Z$ is free of rank $|\cI|$ as a $\Sym U$-module, we call $\cZ$ {\bf free}.
When there is no chance for confusion, we may refer to $Z$ itself as a localization algebra.
\end{definition}


\begin{example}\label{geom-gm}
By the preceding discussion, the equivariant cohomology ring 
$H^*_T(X)$ carries a natural structure of 
a localization algebra, with $U = \mathfrak{t}^*$ and $\cI = X^T$.  
It is both strong and free if and only if $H^{\operatorname{odd}}(X)=0$.
\end{example}

If $H^{\operatorname{odd}}(X)=0$, it is often easier to think of the morphisms  of \eqref{algebras} 
in terms of the dual morphisms of schemes:
$$\t\twoheadleftarrow H^T_2(X)\leftarrow \Spec H_T^*(X)\twoheadleftarrow X^T\times \t.$$
The composite map $\Spec H_T^*(X)\to \t$ is a flat family of schemes,
with zero fiber equal to the fat point $\Spec H^*(X)$ and with general fiber isomorphic to 
$\Spec H^*(X^T)\cong X^T$.
For each $\a\in X^T$, let $H_\a$ be the image of $\{\a\}\times\t$ in $H^T_2(X)$,
a linear subspace that projects isomorphically onto $\t$.  The union of all of these subspaces
is equal to the spectrum of the subring of $H^*_T(X)$ that is generated by the degree two part $H^2_T(X)$; 
equivalently, it is the image of the map from $\Spec H_T^*(X)$ to $H^T_2(X)$ \cite[3.2]{GM}.
This leads us naturally to the following definition, which can be found
in \cite[\S 8.1]{GM}.

\begin{definition}\label{fibered}
A {\bf fibered arrangement} is a surjective map of finite-dimensional complex 
vector spaces $E\twoheadrightarrow F$ along with a finite set $\cI$ and a collection
$\{H_\a\mid\a\in\cI\}$ of linear subspaces of $E$ that project isomorphically onto $F$.
For example, a localization algebra $(U,Z,\cI,h)$ gives rise to a fibered arrangement
by taking $E = Z_2^*$, $F=U^*$, and $H_\a$ equal to the image of the dual of the degree $2$ 
part of the $\a$
component of the localization map $h$.
\end{definition}

\begin{example}\label{ptwo}
  Let $X = \mathbb{P}^2$, and let $T\subset\operatorname{PGL}_2$ be
  the diagonal subgroup.  Then $T$ acts on $X$ with three isolated fixed
  points.  The ring $H_T^*(X)$ is isomorphic to
  $\C[b_1,b_2,b_3]/\langle b_1b_2b_3\rangle$, where $b_i \in H^2_T(X)$ is 
  the equivariant Thom class of the $i^\textrm{th}$ coordinate projective line
  in $X$.  The subring $\Sym\t^*$ is generated by the classes
  $b_1-b_2$ and $b_2-b_3$.  The vector space $H^T_2(X)$ is
  3-dimensional, with coordinates $b_1$, $b_2$, and $b_3$.  The kernel
  of the map to $\t$ is generated by the $(1,1,1)$ vector, and the
  three subspaces $H_\a$ are the coordinate hyperplanes.
\end{example}

\begin{example}\label{two lines}
Let $X$ be a pair of projective lines touching at a single point, and let $T$
be a one-dimensional torus.  We consider the action of $T$ on $X$ such
that $T$ acts effectively on each component, and the double point is an attracting
fixed point for one component and a repelling fixed point for the other.
The ring $H_T^*(X)$ is isomorphic to 
$\C[c_1,c_2,c_3]/\langle c_1c_2, c_1c_3, c_2c_3\rangle$, where $c_i$ is a degree 2
generator whose restriction to the fixed point $p_j$ is $\delta_{ij}$
times a fixed generator of $H^2_T(pt)$.
The vector space $H^T_2(X)$ is 3-dimensional, with coordinates 
$c_1$, $c_2$, and $c_3$.
The kernel of the map to $\t$ is defined by the equation $c_1+c_2+c_3=0$, and the
three subspaces $H_\a$ are the coordinate lines.
\end{example}

Examples \ref{ptwo} and \ref{two lines} motivate the notion of dual fibered
arrangements and dual localization algebras.

\begin{definition}\label{duality}
Consider a fibered arrangement with notation as in Definition \ref{fibered}.
Its {\bf dual} is given by $E^*\twoheadrightarrow E^*/F^*$, along with the 
linear subspaces $H_\a^\perp\subset E^*$, indexed by the same finite set $\cI$.
A {\bf duality} between two localization algebras $\cZ$ and $\cZ^\vee$ is an isomorphism
between the fibered arrangement associated to $\cZ^\vee$ and the dual of the fibered
arrangement associated to $\cZ$.
Thus it consists of an identification of $\cI^\vee$ with $\cI$ and a perfect pairing
between $Z_2^*$ and $(Z_2^\vee)^*$ such that each $H_\a\subset Z_2^*$ is the perpendicular
space to $H_\a^\vee\subset (Z_2^\vee)^*$, and the kernels of the projections to $U^*$
and $(U^\vee)^*$ are also perpendicular to each other.
\end{definition}

\begin{example}\label{first duality}
The localization algebras in Examples \ref{ptwo} and \ref{two lines} are dual via
the perfect pairing of vector spaces with respect to which
$b_1, b_2, b_3$ and $c_1, c_2, c_3$ are dual coordinate systems, and the bijection of
fixed point sets that takes $L_i\cap L_j$ to $p_k$ for $i,j,k$ distinct.
\end{example}

\section{Koszul, quasi-hereditary, and standard Koszul algebras}
\label{Koszul preliminaries}
In this section we review the well-known definitions of quadratic, Koszul, and quasi-hereditary
algebras, along with the slightly less well-known notion of a standard Koszul algebra.
Let $\cI$ be a finite set of order $n$, and
let $R := \C\{e_\a\mid\a\in\cI\}$ be a ring spanned by pairwise
orthogonal idempotents.
Let $M$ be a finitely generated $R$-bimodule,
and let $W \subset M \otimes_R M$ be a sub-bimodule.

Let $$A := T_R(M)\big/\langle W\rangle$$
be the associated quadratic algebra.  For all $\a\in\cI$, let 
$$L_\a := A\Big{/}A_+\oplus\C\{e_\b\mid\b\neq\a\}$$ be the simple right $A$-module indexed by $\a$,
and let $P_\a := e_\a A$ be its projective cover.

\begin{definition}
A complex $\cdots \to M_{i+1}\to M_{i}\to
M_{i-1}\to\cdots$ of graded right $A$-modules is called {\bf linear} if $M_i$ is generated in degree $i$.  
The algebra $A$ is called {\bf Koszul} if each simple module $L_\a$
admits a linear projective resolution.
\end{definition}

Suppose that we are given a partial order $\leq$ on $\cI$, and consider the idempotents
$$\ep_\a := \sum_{\gamma\not\leq\a}e_\gamma\,\,\,\text{and}\,\,\,\ep_\a' := \ep_\a + e_\a.$$  
The {\bf right-standard module}
$V_\a$ is defined to be the largest quotient of $P_\a$ that is supported at or below $\a$,
that is $$V_\a := e_\a A \Big{/}e_\a A \,\ep_\a A.$$
Left-projective modules and left-standard modules are defined similarly.

\begin{definition}
Consider the natural surjections 
$$P_\a\overset{\Pi_\a}{\longrightarrow} V_\a\overset{\pi_\a}{\longrightarrow}L_\a.$$
The algebra $A$ is called {\bf quasi-hereditary} if the following two conditions hold for all $\a\in\cI$:
\begin{itemize}
\item $\ker\pi_\a$ admits a filtration with each subquotient isomorphic to $L_\b$ for some $\b<\a$
\item  $\ker\Pi_\a$ admits a filtration with each subquotient isomorphic to $V_\gamma$ for some $\gamma>\a$.
\end{itemize}
\end{definition}

\begin{remark}
This is equivalent to asking that the regular right $A$-module admits a filtration with standard subquotients, and that the endomorphism algebra of each right-standard module 
$V_\a$ is a division algebra \cite[\S 1]{ADL03}.  It is also equivalent to requiring that the standard modules form an exceptional sequence with respect to the partial order on $\cI$ \cite[Proposition 2]{Bez03}.
\end{remark}

\begin{definition}\label{def:stdk}
The algebra $A$ is called {\bf standard Koszul} if it is finite-dimensional and 
quasi-hereditary,\footnote{\cite{ADL03} is internally inconsistent
about whether or not a standard Koszul algebra should be required to be either
finite-dimensional or quasi-hereditary.  For the purposes of this paper, we require both.}
each right-standard module $V_\a$ admits a linear projective resolution, and the analogous
condition holds for left-standard modules, as well.\end{definition}


\begin{theorem}{\em\cite[Theorem 1]{ADL03}}
If $A$ is standard Koszul, then it is Koszul.
\end{theorem}

The {\bf quadratic dual} $A^!$ of $A$ is defined as the quotient 
$$A^! := T_RM^*/\langle W^\perp\rangle,$$ where $M^* = \Hom_\C(M,\C)$, with the
natural $R$-bimodule structure for which $e_\a M^* e_\b = (e_\b M e_\a)^*$.
It is a well-known fact that $A^!$ is Koszul if and only if $A$ is.
The analogous fact holds for standard Koszulity, as well.

\begin{theorem}{\em\cite[Theorem 3]{ADL03}}
If $A$ is standard Koszul, then $A^!$ is standard Koszul with respect to the opposite
partial order on $\cI$.
\end{theorem}

\begin{remark}
In fact, it is shown in \cite[Theorem 3]{ADL03} that any finite-dimensional,
quasi-hereditary, Koszul algebra is standard Koszul if and only if its dual is quasi-hereditary.
Thus standard Koszul algebras form
the largest class of simultaneously Koszul and quasi-hereditary finite-dimensional algebras 
that is closed under the operation of quadratic duality.
\end{remark}

We conclude with two technical lemmas that we will need in Section \ref{sec:flexible-algebr}.
The first says that if we express a standard Koszul algebra $A$ as 
a quadratic quotient of the path 
algebra of a quiver with vertex set $\cI$, that quiver has no loops, and it only has
arrows between nodes that are comparable in our partial order.  
The second says that any path of length 2 that starts and ends at $\a$
may be uniquely expressed as a sum of paths that avoid all nodes that lie below $\a$.

\begin{lemma}\label{noloops}
If $A$ is standard Koszul and $e_\a M e_\b \neq 0$,
then either $\a <\b$ or $\b < \a$.  In particular, $e_\a M e_\a = 0$ for all $\a\in\cI$.
\end{lemma}

\begin{proof}
For any $A$-module $N$, the {\bf cosocle} of $N$ is defined
to be the largest semisimple quotient of $N$.
Consider the right $A$-module $N = (P_\a)_{\geq 1}/(P_\a)_{\geq 2}$,
which is isomorphic as an $R$-module to $e_\a M$.
Since $N$ has a grading that is concentrated in a single
degree, it is semisimple, and therefore a quotient of the cosocle
of $(P_\a)_{\geq 1}$.

The standard filtration of $P_\a$ induces a filtration of $(P_\a)_{\geq 1}$
with each subquotient isomorphic to either $\ker(\pi_\a)$ or $V_\ga$ for some $\ga>\a$.
This in turn induces a filtration of the cosocle of $(P_\a)_{\geq 1}$, with subquotients
isomorphic to the cosocle of $\ker(\pi_\a)$ or of $V_\ga$ for some $\ga>\a$.
We know that $\ker(\pi_\a)$ only has composition factors of the form $L_\b$ for $\b<\a$,
and that the cosocle of $V_\ga$ is isomorphic to $L_\ga$.  Thus the simple modules
that appear in the cosocle of $(P_\a)_{\geq 1}$ are all of the form $L_\b$ for $\b<\a$ or $\b>\a$.
Since $N = e_\a M$ is a quotient of the cosocle of $(P_\a)_{\geq 1}$, the same is true for $N$.
\end{proof}

\begin{lemma}\label{nu} If $A$ is standard Koszul, then
the projection $e_\a M\ep_\a\otimes_R \ep_\a Me_\a\to e_\a A_2 e_\a$
is an isomorphism for every $\a\in\cI$.
\end{lemma}

\begin{proof}
Since $A$ is standard Koszul, it is {\bf lean} \cite[1.4]{ADL03}, which means that
$$\ep'_\a(\rad A)\ep'_\a(\rad A)\ep'_\a = \ep'_\a(\rad A)^2\ep'_\a$$
for all $\a\in\cI$.  Since $A$ is finite-dimensional by assumption, we have 
$\rad A = A_+$, so multiplying on both the left and the right by $e_\a$ and looking in degree 2,
we have a surjection 
\begin{equation}\label{a projection}
e_\a M\ep_\a\otimes_R \ep_\a Me_\a = e_\a M\ep'_\a\otimes_R \ep'_\a Me_\a\twoheadrightarrow
e_\a A_1 \ep'_\a A_1 e_\a = e_\a A_2 e_\a,
\end{equation}
where the equality on the left follows from Lemma \ref{noloops}.

Consider the exact sequence
\[0 \to K \to P_\a \to V_\a \to 0.\]
Since $V_\a e_\a$ vanishes in positive degrees, we have
$K_2 e_\a = (P_\a)_2 e_a = e_\a A_2 e_\a$.
Note that $K$ has a
\emph{graded} filtration whose subquotients are isomorphic 
to $V_\b(i)$ for $\b > \a$ and $i > 0$, where $(i)$ denotes shifting the
grading up by $i$ --- see \cite[(1.2)]{CPS}.  
It follows that $K_1 \cong \bigoplus_{\b > \a} e_\a M e_\b$.

The subquotients $V_\gamma(j)$ of $K$ with $j \ge 2$ cannot 
contribute to $K_2 e_\a$, so we have an isomorphism
\[K_2 e_\a \cong \bigoplus_{\b > \a} e_\a M e_\b \otimes (V_\b)^{}_1 e_\a = 
\bigoplus_{\b > \a}  e_\a M e_\b \otimes e_\b M e_\a.\]
It follows that $\dim e_\a A_2 e_\a = \dim  e_\a M\ep_\a\otimes_R \ep_\a Me_\a$,
so the projection \eqref{a projection} is an isomorphism.
\end{proof}


\section{Flat deformations of Koszul algebras}
\label{sec:canon-deform-quadr}
In this section we study graded deformations of Koszul algebras.  To begin
we consider a quadratic algebra $A = T_R(M)/\langle W \rangle$
which is not necessarily Koszul.
Let $U$ be a finite-dimensional $\C$-vector space and let
$S = \Sym U$, graded so that elements of $U$ have degree two.  
Now suppose given a graded deformation of $A$ over $U^*$, that is,  
a graded $R$-algebra $\tilde{A}$ together with
graded homomorphisms
\[S \stackrel{j}{\lra} \tilde{A} \stackrel{\pi}{\lra} A\]
so that $j$ maps into the center of $\tilde{A}$ and $\pi$ 
induces an isomorphism $A \cong \tilde{A}/\langle j(U)\rangle$
of graded algebras. 

Since $S$ is generated in degree two,
the map $\pi$ is an isomorphism in degrees zero and 
one, so we have $\tilde{A}_0 \cong R$ and $\tilde{A}_1 \cong M$.
In degree two we have the right exact sequence
\[ R \otimes U  \to \tilde{A}_2 \to A_2 \to 0.\]
We will make the additional assumption that $\tilde{A}$ is 
flat over $S$; in particular, this implies that the 
sequence above is in fact short exact. 

Since $W$ is contained in the kernel of the projection
$T_R(M) \to A$, we get a map of $R$-bimodules
$\Phi\colon W \to R \otimes U$ by letting $\Phi(w)$ be the
image of $w$ under the multiplication map
$T_R^2(M) = T_R^2(\tilde{A}_1) \to \tilde{A}$.  

Note that if $N$ is any $R$-bimodule and $V$ is a $\C$-vector
space, then any map $\Psi\colon N \to R \otimes V$ of $R$-bimodules
must annihilate any ``off-diagonal" summand
$e_\a N e_\b \subset N$ with $\a \ne \b$.  In particular,
this means that such bimodule maps are in bijection 
with linear maps $\Psi^\circ\colon N \to V $ that
kill the off-diagonal summands, via the formula
\[\Psi(e_\a x e_\b) = e_\a \otimes \Psi^\circ(e_\a x e_\b).\]


For any map $\Psi\colon W \to R \otimes U$ of $R$-bimodules, we define
\begin{equation}\label{psi-def}
\tilde{A}_\Psi := T_R(M) \otimes S \big/\langle w \otimes 1 - 1\otimes \Psi(w) 
\mid w \in W\rangle.
\end{equation}
We then get a surjective map $\tilde{A}_\Phi \to \tilde{A}$ of graded algebras
which becomes an isomorphism upon tensoring over $S$ with $\C$.  Since
we are assuming that $\tilde{A}$ is flat over $S$, this map must be  
an isomorphism even before tensoring.  
Thus every flat graded deformation of $A$ with deformation parameters in degree two
arises from a bimodule map $\Psi\colon W\to R \otimes U$
as in Equation \eqref{psi-def}.  It is not the case, however, that every bimodule
map gives a flat graded deformation.  If $A$ is Koszul, we have the following criterion for flatness.

\begin{theorem}\label{flat}  Suppose that $A$ is Koszul.
Let $\Psi: W\to R \otimes U$ be a bimodule map, and let $\tilde{A}_\Psi$ be the graded deformation
of $A$ given in Equation \eqref{psi-def}.  This deformation is flat if and only if
$\Psi^\circ$ factors through the quotient map $W \cong (A^!_2)^* \to Z(A^!)_2^*$.
\end{theorem}

\begin{remark}\label{universal remark}
Since $Z(A^!)_2$ has no off-diagonal summands,
Theorem \ref{flat} implies that graded flat deformations of $A$ over $U^*$ 
with deformation parameters in degree two
are in bijection with linear maps $\psi\colon Z(A^!)_2^* \to U$.
If we take $U$ to be $Z(A^!)_2^*$ and $\psi$ to be the identity map,
we call the resulting ring $\tilde{A}$
the {\bf universal deformation} of $A$.
It is universal in the sense that if 
$\psi' \colon Z(A^!)_2^* \to U'$ is another linear map,
then the corresponding deformation is isomorphic to
$\tilde{A} \otimes_S \Sym U'$, where the map
$S \to \Sym U'$ is induced by $\psi'$.

In Section \ref{sec:deformed O}, we will also need the analogous
statement for deformations over power series rings.
If $\wh A$ is a flat deformation of $A$ over the spectrum of the power
series ring $\prod_{i=0}^\infty \Sym^i U'$ and $\wh A$
admits a {\bf formal grading} $\wh A = \prod_{i=0}^\infty \wh A_i$ compatible with its algebra structure, then $\wh A$ 
may be obtained from the universal deformation of $A$ via a base change 
to $U^*$ followed by a completion at the unique graded maximal ideal of $\Sym U$.
\end{remark}

\begin{proof}[Proof of Theorem \ref{flat}]
For any linear map $\chi\colon U \to \C$ 
let $\C_\chi$ denote the associated 
one-dimensional $\Sym U$-algebra, and
consider the specialization 
\[A_\chi := \tilde A_\Psi \otimes_{\Sym U} \C_\chi\]
of $\tilde A_\Psi$ at the point $\chi \in U^*$.
Explicitly, the ring $A_\chi$ 
is the quotient of the tensor algebra $T_RM$ by 
the two-sided ideal 
$\langle w - \chi^{}_R \circ \Psi(w) \mid w \in W\rangle$,
where $\chi_R\colon R\otimes U \to R$ is  
given by $\chi_R(e_\a \otimes u) = e_\a\chi(u)$ for all $\a\in\cI$ and $u\in U$.
The grading on $\tilde A_\Psi$ induces a filtration
on $A_\chi$, and the ring $\tilde A_\Psi$ is
flat if and only if the natural surjection
$A \to \gr A_\chi$ is an isomorphism for
all $\chi$.

Our theorem now follows directly from a result of Braverman and Gaitsgory \cite{BG96}.
They study a more general situation, taking the quotient $Q$ of
the tensor algebra $T_R(M)$ by the two-sided ideal 
$\langle w - a(w) - b(w) \mid w \in W\rangle$, 
where $a \colon W \to M$
and $b\colon W \to R$ are maps of $R$-bimodules.  
Their main result\footnote{Braverman and Gaitsgory only treat the case 
when $R$ is a field, but their results easily generalize to our semisimple ring $R$.  
Their condition (I) forces $a$ and $b$ to be bimodule
maps.} \cite[4.1]{BG96} gives 
necessary and sufficient conditions on $a$ and $b$ 
to have $\operatorname{gr} Q \cong A$.  In our case,
we have $a = 0$ and $b = \chi_R \circ \Psi$.  
In this situation their conditions reduce to the statement that the map
\[b \otimes \id - \id \otimes\, b \colon (W \otimes_R M) \cap (M \otimes_R W) \to M\]
vanishes.

To relate this condition to the dual ring,  
note that 
\[A^!_3 = M^* \otimes_R M^* \otimes_R M^* \big/ (W^\perp \otimes_R M + M \otimes_R W^\perp)\]  
is naturally dual to 
$(W \otimes_R M) \cap (M \otimes_R W)$.  Thus we may identify $Z(A^!)_2$ with the set of 
$\C$-linear maps $b^\circ\colon W \to \C$ that kill the off-diagonal terms 
(this implies that $b^\circ$ commutes with the idempotents) and for which 
\[b^\circ \otimes \gamma - \gamma \otimes b^\circ \colon (W \otimes_R M) \cap (M \otimes_R W) \to \C\]
vanishes for any $\gamma \in M^*$.  Here we consider $W \otimes_R M$ as 
a subspace of $W \otimes M$ in the obvious way, so 
$(b^\circ \otimes \gamma)(w \otimes m) = \sum_\a b^\circ(e_\a w e_\a)\gamma(e_\a m)$, and 
similarly for $\gamma \otimes b^\circ$. 
Then if $b\colon W \to R \otimes \C = R$ is the bimodule map corresponding
to a $\C$-linear map $b^\circ$, we have  
$b^\circ \otimes \gamma - \gamma \otimes b^\circ = \gamma \circ (b \otimes \id - \id \otimes b)$, 
so $b^\circ$ represents an element of $Z(A^!)_2$ if and only if $b \otimes \id - \id \otimes b = 0$.  

In other words, we have shown that $\tilde{A}_\Psi$ is flat if and only if 
$(\chi_R \circ \Psi)^\circ = \chi \circ \Psi^\circ$ is
central for all $\chi\in U^*$, which is equivalent to saying that $\Psi$ factors through
$Z(A^!)^*_2$.
\end{proof}



\begin{remark}
If we drop the Koszulity hypothesis, then the ``if" direction of 
Theorem \ref{flat} becomes false; it fails, for example, when $|\cI| = 1$ and 
$A = \C\langle x,y\rangle / \langle x^2, y^2 - xy\rangle$.\footnote{We thank
Andrew Connor for this example.}
\end{remark}

\begin{remark}

Theorem \ref{flat} can be understood more abstractly using the fact that 
Koszul duality induces an equivalence of derived categories of graded modules  
$$D^b(A-\gmod)\cong D^b(A^!-\gmod).$$  Since Hochschild cohomology is equal to the Ext-algebra 
of the identity functor on the derived category \cite[1.6]{Toen}, this also induces an isomorphism of Hochschild cohomology groups.  
The behavior of this equivalence on grading shift functors is such that
the group $H\! H^r(A^!)_{s}$ is identified with $H\! H^{r + s}(A)_{-s}$
\cite[1.2.6]{BGS96}.  
In particular,
if $A$ is Koszul, $Z(A^!)_2 = H\! H^0(A^!)_2$ is naturally isomorphic to $H\! H^2(A)_{-2}$.
In fact, the proof of the main result of \cite{BG96} proceeds by showing that
an element of $Z(A^!)_2$ lifts to 
an $R$-bimodule map $M\otimes M \to R$ satisfying a cocycle
condition which allows it to represent a class in $H\! H^2(A)_{-2}$.
\end{remark}

\section{Deformed standard modules and malleable algebras}\label{sec-malleable}
Throughout this section we assume that $A$ is a standard Koszul algebra,
$S = \Sym U$ is a polynomial ring, and $\tilde A$ is a 
flat graded deformation of $A$ over $U^*$.
Consider the right $\tilde A$-modules
$$\tilde P_\a := e_\a\tilde A\,\,\,\,\,\,\text{and}\,\,\,\,\,\,\tilde V_\a := \tilde P_\a 
\Big{/}e_\a \tilde A \,\ep_\a \tilde A.$$  Since $\tilde P_\a$ is a summand of $\tilde A$,
it is a flat deformation of $P_\a$.
The purpose of this section is to show that the center of $\tilde A$ acts on each $\tilde V_\a$
via a central character
$$
	h_\a: Z(\tilde A) \longrightarrow S.
$$  
It will follow that the data  $\cZ(\tilde A):= (U,\, Z(\tilde A),\, \cI,\, \oplus h_\a)$
form a localization algebra.

For any $\a \in \cI$, consider the algebra $C_\a :=\ep_\a'A\,\ep_\a'\subset A$ and its deformation
$\tilde C_\a := \ep_\a'\tilde A\,\ep_\a'\subset\tilde A$.  By \cite[3.9]{ADL03}, $C_\a$ is
standard Koszul.  The deformed algebra $\tilde C_\a$ is a direct summand of $\tilde A$ as an $S$-module,
so it is flat over $U^*$.
For any $\a\leq\b$, let
$$V^\a_\b := e_\b C_\a\Big{/}e_\b C_\a \ep_\b C_\a$$ be the standard cover of 
$L_\b$ in the category of right $C_\a$-modules, and consider its deformation
$$\tilde V^\a_\b := e_\b\tilde C_\a\Big{/}e_\b\tilde C_\a\ep_\b\tilde C_\a.$$

\begin{lemma}\label{st-induced}
We have an isomorphism $V^\a_\b\otimes_{C_\a}\ep_\a'A\cong V_\b$
of right $A$-modules, and an isomorphism 
$\tilde V^\a_\b\otimes_{\tilde C_\a}\ep_\a'\tilde A\cong \tilde V_\b$ of right $\tilde A$-modules.
\end{lemma}

\begin{proof}
Using the equalities $e_\b\ep_\a' = e_\b$ and $\ep_\b\ep_\a' = \ep_\b = \ep_\a' \ep_\b$,
we have 
\begin{eqnarray*}
V^\a_\b\otimes_{C_\a}\ep_\a'A &=& 
\left(e_\b C_\a\Big{/}e_\b C_\a \ep_\b C_\a\right)\otimes_{C_\a}\ep_\a'A\\
&=& \left(e_\b A\,\ep_\a'\Big{/}e_\b A\ep_\b A\,\ep_\a'\right)
\otimes_{C_\a}\ep_\a'A\\
&\cong& e_\b A\Big{/}e_\b A\,\ep_\b A = V_\b.
\end{eqnarray*}
The proof of the second statement is identical.
\end{proof}

\begin{remark}
The most important case of Lemma \ref{st-induced}, and also the easiest one to think about,
is the case in which $\a=\b$, so that $V_\b^\a = V_\a^\a$ is isomorphic to the simple module
for $C_\a$ supported at the node $\a$.
\end{remark}

\begin{proposition}\label{standard def}
For all $\a\in\cI$, $\tilde V_\a$ is a flat deformation of $V_\a$.
\end{proposition}

\begin{proof}
We first consider the case where $\a$ is a minimal element of our poset.  
In this case $V_\a \cong L_\a$
is one-dimensional.  Then by semicontinuity, 
the $\tilde A_{\chi}$-module 
$(\tilde V_\a)_{\chi} = e_\a\tilde A_{\chi}\Big{/}e_\a \tilde A_{\chi} \,\ep_\a \tilde A_{\chi}$
has dimension 0 or 1 for every $\chi \in U^*$.  We must show that that dimension
is equal to 1 for every $\chi$, or equivalently that $e_\a\notin e_\a\tilde A_{\chi} \,\ep_\a \tilde A_{\chi}e_\a$.
Since $\a$ is minimal, Lemma \ref{nu} tells us that there are no nontrivial relations
among loops of length 2 based at $\a$.  In particular there are no relations to deform,
and the conclusion follows.

In the general case, $\a$ is a minimal element of the poset of simples for the subalgebra
$C_\a$, so $\tilde V^\a_\a$ is a flat deformation of $V^\a_\a$.
This then implies the result for $\tilde V_\a \cong \tilde V_\a^\a\otimes_{\tilde C_\a}\ep'_\a\tilde A$.
\end{proof}

\begin{lemma}\label{projective def}
The regular right $\tilde A$-module $\tilde A$
admits a filtration indexed by $\cI$ for which
the subquotient indexed by $\a \in \cI$ is
isomorphic to a direct sum of shifts of
the deformed standard $\tilde V_\a$.
\end{lemma}

\begin{proof}
Choose a maximal index $\a\in\cI$.  Let $B=A/A e_\a A$,
and let $W_B\hookrightarrow B_1\otimes_R B_1$
be the space of relations for $B$.  
The algebra $B$ is standard Koszul by \cite[3.9]{ADL03}.
Let $\tilde B = \tilde A/\tilde Ae_\a\tilde A$.  We
first use the results of Section \ref{sec:canon-deform-quadr}
to show that $\tilde B$ is a flat deformation of $B$ over $U^*$.

Let $\bar{e} = 1 - e_\a$.
The surjection $A \surj B$ induces an inclusion $B^! \hookrightarrow A^!$
with image contained in $\bar{e}A^!\bar{e} \subset A^!$.
By \cite[2.5]{ADL03}, the image is equal to $\bar{e}A^!\bar{e}$.
The map $q \colon A^!_2 \to B^!_2$ given by $q(x) = \bar{e}x\bar{e}$ 
is a left inverse for the inclusion $B^!_2 \hookrightarrow A^!_2$.    
It follows that if the deformation $\tilde{A}$ is described by an
$A_0$-bimodule map $\Psi\colon W_{\! A} \to A_0 \otimes U$ as in Section 
\ref{sec:canon-deform-quadr}, then $\tilde{B}(\cV)$ is 
described by the composition
\[W_{B} \stackrel{q^*}{\lra} W_{\! A} \stackrel{\Psi}{\lra} A_0 \otimes U \stackrel{\cdot \bar{e}}{\lra} B_0 \otimes U.\]
Then $q$ sends $Z(A^!)_2$ into $Z(B^!)_2$,
so Theorem \ref{flat} implies that $\tilde{B}(\cV)$ is flat over $U^*$.

A deformed standard module over $\tilde B$ becomes a deformed standard $\tilde A$-module
under the quotient homomorphism $\tilde A \to \tilde B$.
Thus, inducting on the size of our poset, we may assume that the right $\tilde A$-module
$\tilde B$ has a filtration by graded $\tilde A$-modules 
indexed by $\cI \setminus \{\a\}$
such that the associated graded indexed by $\b$ 
is isomorphic to a direct sum of shifts of $\tilde V_\b$.

Consider the exact sequence
\begin{equation}\label{induct}
0\to\tilde A e_\a \tilde A \to \tilde A \to \tilde B\to 0.
\end{equation}  
We have $$\tilde A e_\a \tilde A=\tilde A e_\a\otimes_{e_\a
    \tilde A e_\a} e_\a \tilde A=\tilde A e_\a\otimes_{e_\a
    \tilde A e_\a} \tilde P_\a=\tilde A e_\a\otimes_{e_\a
    \tilde A e_\a} \tilde V_\a,
$$
where the last equality follows from the maximality of $\a$.
Maximality also implies that
$e_\a \tilde A e_\a\cong S$ (by Lemma \ref{nu}),
thus the $\tilde A$-module $\tilde A e_\a \tilde A$ is isomorphic to 
a direct sum of $\dim Ae_\a$ copies of $\tilde V_\a$.  The result 
then follows from the exact sequence \eqref{induct}.
\end{proof}

\begin{corollary}\label{deformed standard characterization}
Suppose that $\tilde V$ is a graded $\tilde{A}$-module
which is a flat deformation of $V_\a$; that is, it is free as an
$S$-module and $\tilde{V} \otimes_S \C_0$ is isomorphic to $V_\a$ as 
an $A$-module.  Then $\tilde{V} \cong \tilde{V}_\a$.
\end{corollary}

\begin{proof}
Consider the surjection
$\tilde{V} \surj \tilde{V} \otimes_S \C_0 \cong V_\a$.
Since $\tilde{P}_\a$ is projective, we can lift
the map $\tilde{P}_\a \to \tilde{V}_\a \to V_\a$ to
a map $\phi\colon\tilde{P}_\a \to \tilde{V}$.  The 
fiber of $\phi$ over $0$ is a surjection, so Nakayama's lemma tells us that
$\phi$ itself must be surjective.
For any $\b \in \cI$, the natural map
$$\Hom_{\tilde{A}}(\tilde{P_\b}, \tilde{V}) \otimes_S \C_0 \to \Hom_A(P_\b, \tilde{V}\otimes_S \C_0)$$
is injective, thus $\Hom_{\tilde{A}}(\tilde{P_\b}, \tilde{V}) = 0$ for any $\beta \not\le \alpha$.
It follows that
$\ker \phi$ contains $e_\a\tilde{A}\,\ep_\a\tilde{A}$, 
hence $\tilde{V}$ is a quotient of $\tilde{V}_\a$.  
Since they are both free $S$-modules and 
their fibers over $0$ are isomorphic, we have 
$\tilde{V} \cong \tilde{V}_\a$.
\end{proof}

The next result says that our deformed standard modules have well-defined 
central characters. 
\begin{proposition}\label{hadef}
  For each $\a\in\cI$, there is an $S$-algebra homomorphism $h_\a:Z(\tilde A)\to S$ such
  that $v\cdot \z=h_\a(\z)v$ for all $v\in \tilde V_\a$ and $\z\in Z(\tilde A)$.
  \end{proposition}

\begin{proof}
If $\a$ is minimal, then $V_\a\cong L_\a$ and $\tilde V_\a\cong S$ by Proposition
\ref{standard def}, and the claim follows from the fact that $\tilde A$ acts on $\tilde V_\a$
by $S$-module endomorphisms.  The general case follows from
the isomorphism $\tilde V_\a\cong\tilde V_\a^\a\otimes_{\tilde C_\a}\ep_\a'\tilde A$ of Lemma
\ref{st-induced}.
\end{proof}


As a consequence, we can construct a localization algebra from a flat graded deformation 
$\tilde{A}$ of a standard Koszul algebra $A$, thus completing the proof of Theorem \ref{first}.

\begin{corollary}\label{gma}
The data $\cZ(\tilde A):= (U,\, Z(\tilde A),\, \cI,\, \oplus h_\a)$
form a localization algebra.
\end{corollary}

Let $K$ be the fraction field of $S$, and for any $S$-module $N$,
let $N^\infty = N\otimes_S K$.  The filtration of $\tilde A$
from Lemma \ref{projective def}
induces a filtration of $\tilde A^\infty$ with subquotients isomorphic to direct sums
of modules of the form $\tilde V_\a^\infty$.

\begin{theorem}\label{remarkable-thm}  
Suppose that $A\cong A^{op}$ as $R$-algebras.
  The following are equivalent:
  \begin{enumerate}
  \item The deformed standard filtration of $\tilde A^\infty$ splits.
  \item 
The action map $\tilde A^\infty\to\bigoplus_\a\operatorname{End}_K(\tilde V_\a^\infty)$ is an isomorphism.
  \item The map $\bigoplus_\a h^\infty_\a:Z(\tilde A)^\infty\to \bigoplus_{\a} K$ 
  is an isomorphism.
 \item The maps $\{h_\a\mid\a\in\cI\}$ from $Z(\tilde A)$ to $S$ are all distinct.  
\end{enumerate}
\end{theorem}

\begin{proof}
 $(1)\Rightarrow(2)$:
  Since the deformed standard filtration of $\tilde A^\infty$ splits, 
  an element of $\tilde A^\infty$ which kills every $\tilde V_\a^\infty$ must also act trivially on $\tilde A^\infty$, and must therefore be zero.  Thus the action map is injective.  
Comparing dimensions, we have 
  \begin{equation*}
    \dim_K\tilde A^\infty=\dim_\C A=\sum_{\a,\b}[L_\a: P_\b]=\sum_{\a,\b,\ga}[L_\a:V_\ga][V_\ga:P_\b]=\sum_{\a,\b,\ga}[L_\a:V_\ga][L_\b:V_\ga]=\sum_\ga (\dim V_\ga)^2,
  \end{equation*}
  where the 
  penultimate equality follows from BGG reciprocity \cite[3.11]{CPS88} and the isomorphism $A\cong A^{\text{op}}$.  Thus our map must be an isomorphism.

$(2)\Rightarrow (3)$: This follows from the fact that $Z(\tilde A)^\infty \cong Z(\tilde A^\infty)$.  

$(3)\Rightarrow (4)$: Immediate from surjectivity.

$(4)\Rightarrow (1)$: If all of the deformed standard modules have different central
characters, then the filtration of $\tilde A$ can be split by
taking the isotypic decomposition for the action of the center.
Here we use the fact that there is no overlap between central characters of consecutive
subquotients, which follows from Lemma \ref{projective def}.
\end{proof}

\begin{definition}\label{malleable}
If $\tilde{A}$ satisfies the conditions of Theorem \ref{remarkable-thm}, 
we will call $\tilde A$ {\bf malleable}.
\end{definition}

\begin{proposition}\label{malleable implies strong}
Suppose that $A\cong A^\text{op}$.  Then
$\cZ(\tilde A)$ is a strong localization algebra
iff $\tilde A$ is malleable.
\end{proposition}

\begin{proof}
By Definition \ref{gmalg}, $\cZ(\tilde A)$ is strong if and only if the kernel and cokernel of $\oplus_\a h_\a$
are torsion.  This is equivalent to asking that $\oplus_\a h_\a^{\infty}$
be an isomorphism, which is condition (3) above.
\end{proof}

\section{Flexible algebras}
\label{sec:flexible-algebr}
In this section we define and study flexible algebras in preparation for the next section, 
which contains proofs of our central results, Theorem \ref{bang-flexible}  and Corollary \ref{main}.
Let $A$ be a standard Koszul algebra, let $S=\Sym U$ be a polynomial ring generated in degree 2,
and let $\tilde A$ be a flat graded 
deformation of $A$ over $U^*$.

\begin{definition}
We say that $\tilde A$ is {\bf flexible} if 
the natural projection $Z(\tilde A)_2\to Z(A)_2$ is surjective. 
\end{definition}

\begin{example}\label{two nodes}
Let $\cI = \{1, 2\}$ be the nodes of a quiver with
$r>0$ arrows $x_1,\ldots,x_r$ from 1 to 2 and $s>0$ arrows
$y_1,\ldots,y_s$ from 2 to 1.  Let $A_{rs}$ be the path algebra
modulo the quadratic 
relations $y_jx_i = 0$ for all $i$ and $j$.
Thus a right $A_{rs}$-module is a representation of the quiver for which every
loop based at the node 2 acts trivially.
We have standard modules $V_1 = L_1$ and $V_2 = P_2$ with $\ker\Pi_1\cong P_2^{\oplus r}$
and $\ker \pi_2 \cong L_1^{\oplus s}$, so $A$ is quasi-hereditary with respect
to the order $1<2$.  It is clear that both standard modules have linear projective resolutions.
Furthermore, since the opposite algebra of $A_{rs}$ is isomorphic to $A_{sr}$, 
the same is true for
left-standard modules.  Hence $A_{rs}$ is standard Koszul.  It is also easy to see that
$A_{rs}$ is isomorphic to its own quadratic dual.

The center $Z(A_{rs})$ is spanned by the unit and the $r\times s$ elements $x_iy_j$.
The universal deformation $\tilde{A}_{rs}$ has central generators $u_{ij}$ in degree 2
and relations $y_jx_i = u_{ij}e_2$.
It is easy to check that the generators $x_iy_j$ of $Z(A_{rs})_2$ lift to central elements
of $\tilde A_{rs}$ if and only if $r=s=1$, thus only $\tilde A_{11}$ is flexible.
We note that $\tilde A_{11}$ is also malleable in the sense of Definition \ref{malleable}.
\end{example}

For the rest of this section, let us suppose that our algebra 
$A$ is \textbf{connected}, meaning that
for all $\a, \b \in \cI$ we have $e_\a T_R(M)e_\b\neq 0$.  In other words,
the category of $A$-modules does not split into smaller blocks.  Then
$Z(A)_0 \cong \C$, and so we have an exact sequence 
\begin{equation}\label{exact}
0\to U \to Z(\tilde A)_2 \to Z(A)_2
\end{equation}
which is also right exact if $\tilde A$ is flexible.
In that case, each homomorphism
$h_\a: Z(\tilde A)_2 \to S_2 = U$ of Proposition \ref{hadef} splits this exact sequence.
The difference between any two splittings vanishes on $U$, 
and thus induces a map 
$j_{\a\b}:Z(A)_2\to U$
given by 
\begin{equation*}
j_{\a\b}(z) := h_\b(\z) - h_\a(\z)\,\,\,\,\text{for any lift $\z$ of $z$.}
\end{equation*}

Define maps
\begin{equation}\label{munu}
\mu:Z(A)_2\to \bigoplus_{\a\in\cI}e_\a M\ep_\a\otimes_R \ep_\a Me_\a
\,\,\,\,\,\,\,\,\text{and}\,\,\,\,\,\,\,\,\nu: Z(A)_2\to \tilde A_2
\end{equation}
by setting $\mu(z)$ equal to the unique expression for $z$ as a sum of loops that go first
up and then down (which exists by Lemmas \ref{noloops} and \ref{nu}),
and $\nu(z)$ to the image of $\mu(z)$ in $\tilde A_2$.

\begin{proposition}\label{jabprop}
Suppose that $\tilde A$ is flexible.
For all $z\in Z(A)_2$ and $\tilde a\in e_\a \tilde A e_\b$, we have
$$[\nu(z), \tilde a] = j_{\a\b}(z)\, \tilde a.$$
\end{proposition}

\begin{proof}
Let $z$ be given, and let $\z$ be a lift of $z$ to $Z(\tilde A)_2$.
Since the kernel of the projection from $\tilde A_2$ to $A_2$ is equal to $R\otimes U$,
there exist elements $u_\a\in U$ such that $\z = \nu(z) + \sum_\a u_\a e_\a$.
Since the deformed standard module $\tilde V_\a$ is supported on and below $\a$,
and $\nu(z)$ is expressed in terms of paths that avoid such nodes,
we have $h_\a(\z) = u_\a$ for all $\a\in\cI$, and therefore $j_{\a\b}(z) = u_\b - u_\a$
for all $\a,\b\in\cI$.  Since $\z$ is central, we have
$$[\nu(z), \,\tilde a] = [\nu(z) - \z, \,\tilde a]
= (u_\b - u_\a) \,\tilde a = j_{\a\b}(z) \,\tilde a$$
for all $\tilde a\in e_\a \tilde A e_\b$.
\end{proof}

Proposition \ref{bajprop} may be regarded as a converse to Proposition \ref{jabprop}.

\begin{proposition}\label{bajprop}
Suppose that there exists a collection of linear maps
$$\big\{j'_{\a\b}:Z(A)_2\to U\,\,\big{|}\,\, \a,\b\in\cI\big\}$$
satisfying the following two conditions:
\begin{itemize}
\item $j'_{\a\b} + j'_{\b\gamma} = j'_{\a\gamma}$ for all $\a,\b,\gamma\in\cI$ 
(in particular $j'_{\a\b} = -j'_{\b\a}$ for all $\a,\b\in\cI$)
\item for all $z\in Z(A)_2$ and $\tilde a\in e_\a \tilde A e_\b$, we have
$[\nu(z), \tilde a] = j'_{\a\b}(z)\, \tilde a$.
\end{itemize}
Then $\tilde A$ is flexible, and $j'_{\a\b} = j_{\a\b}$ for all $\a,\b\in\cI$.
\end{proposition}

\begin{proof}
Let $z\in Z(A)_2$ be given; we must show that we can lift it to $Z(\tilde A)_2$.
Choose an element $\delta\in\cI$ arbitrarily, and let 
$$\z = \nu(z) + \sum_{\gamma\in\cI}\, j'_{\delta\gamma}(z) e_\gamma\in \tilde A_2.$$
Then for any $\tilde a\in e_\a \tilde A e_\b$, we have
$$[\z, \tilde a] = [\nu(z), \,\tilde a] + j'_{\delta\a}(z) \,\tilde a - j'_{\delta\b}(z) \,\tilde a
= \Big(j'_{\a\b}(z) +  j'_{\delta\a}(z) - j'_{\delta\b}(z)\Big) \tilde a = 0,$$
where the vanishing of the expression inside the parentheses follows
from the first condition on the homomorphisms $j'_{\a\b}$.
Thus $\z$ is central.

By the same argument that we used in the proof of Proposition
\ref{jabprop}, we have $h_\gamma(\zeta) = j'_{\delta\gamma}(z)$ for all $\gamma$,
thus $j_{\a\b}(z) = h_\b(\zeta) - h_\a(\zeta) = j'_{\delta\b}(z) - j'_{\delta\a}(z) = j'_{\a\b}(z)$.
\end{proof}

We conclude with a lemma that we will use in 
the last section to show that a certain 
flexible deformation of a standard Koszul 
algebra is universal, or at least has the universal 
deformation as a quotient.
Let $\psi\colon Z(A^!)_2^* \to U$ be given, and let
$\tilde A$ be the graded flat deformation provided by 
Theorem \ref{flat}.
Suppose that this deformation is flexible, so that we can define the
maps $j_{\a\b}\colon Z(A)_2\to U$.  

\begin{lemma} \label{image of psi}
For any $\a, \b \in \cI$,
$\Im j_{\a\b}$ is contained in the image of $\psi$.  
\end{lemma}

\begin{proof} 
First suppose that $e_\a M e_\b \ne 0$, and let 
$\tilde{a}$ be any nonzero element of this space.
Then for any $z \in Z(A)_2$,
Proposition \ref{jabprop} implies that $j_{\a\b}(z)\tilde{a}$ 
is a nonzero element of the subring of $\tilde{A}$ generated by
the degree $0$ and $1$ parts, thus $j_{\a\b}(z)$ is in the image of $\Im(\psi)$.
The general case follows from the identity $j_{\a\b} + j_{\b\ga} = j_{\a\ga}$,
since we are assuming that $A$ is connected.
\end{proof}

\section{Koszul duality and GM duality}
\label{sec:Koszul and GM}

In this section we explore the relation between Koszul duality and
the localization algebras of flexible deformations.  
Let $A$ be a connected standard Koszul algebra, and let
$\tilde A$ be its universal 
deformation over $U^* = Z(A^!)_2$.  
On the dual side, let $\tilde A^!$ be the 
universal deformation of the dual ring $A^!$ over 
$(U^!)^* = Z(A)_2$.  Let $S = \Sym U$ and $S^! = \Sym U^!$. 

We will call a standard Koszul algebra flexible if its 
universal deformation is flexible.  If $A$ is flexible, then
we have the maps $j_{\a\b}\colon Z(A)_2 \to U = Z(A^!)_2^*$
constructed in the previous section.  If $A^!$ is flexible,
then the same construction gives maps 
$j^!_{\a\b}\colon Z(A^!)_2 \to U^! = Z(A)_2^*$.

\begin{theorem}\label{bang-flexible}
If $A$ is flexible, then so is $A^!$.  Furthermore, for all $\a,\b\in\cI$,
we have an identity $j^!_{\a\b} = j_{\b\a}^*$ of maps from
$Z(A^!)_2 = U^*$ to $U^! = Z(A)_2^*$.
\end{theorem}

\begin{proof}
By Proposition \ref{bajprop}, it is enough to show that
$$[\nu^!(z^!), \tilde a^!] = j_{\b\a}^*(z^!)\,\tilde a^!$$ for all $z^!\in Z(A^!)_2$
and $\tilde a^!\in e_\a\tilde A e_\b$.  Retracing the computation in the proof
of Proposition \ref{bajprop}, this is equivalent to showing that, for a fixed $\delta\in\cI$,
$$\z^! := \nu^!(z^!) + \sum_{\gamma\in\cI}\, j^*_{\gamma\delta}(z^!) e_\gamma\in \tilde A^!_2$$ 
is central.  It clearly commutes with the idempotents, so it is enough to show that it 
commutes with elements of $M^*$.

We need to work with explicit representatives in the deformed tensor algebra
\[\wt{T_R(M^*)} := T_R(M^*) \otimes \Sym(U^!).\]
Since elements of $U^!$ have degree $2$, we have
\[\wt{T_R(M^*)}_2 = T_R(M^*)_2 \oplus (R\otimes U^!)\,\,\,\text{and}\,\,\, \wt{T_R(M^*)}_3 = T_R(M^*)_3 \oplus (M\otimes U^!).\]
Define a lift $\eta^! \in \wt{T_R(M^*)}$ of $\zeta^!$ by 
 $$\eta^! := \mu^!(z^!) + \sum_{\gamma\in\cI}\, 
j^*_{\gamma\delta}(z^!) e_\gamma \,\in\, T_R(M^*)_2\oplus (R\otimes U^!),$$
where $\mu^!$ is defined as in Equation \eqref{munu}.
Fix a pair of indices $\a,\b\in\cI$, and let $x^!$
be any element of $e_\b M^* e_\a$.  We need to show that the commutator
$[\z^!, x^!]\in T_R(M^*)_3\oplus (M^*\otimes U^!)$ reduces to zero in $\tilde A^!_3$.

Let $Q_1 = M\otimes W$ and $Q_2 = W \otimes M\subset T_R(M)_3$.
By the definition of the quadratic dual, the kernel of the map from
$T_R(M^*)_3$ to $A^!_3$ is equal to
$$M^*\otimes W^\perp + W^\perp \otimes M^* = Q_2^\perp + Q_1^\perp = (Q_1\cap Q_2)^\perp.$$
We now need a similar expression for the relations in $\tilde A^!_3$.
Let $$O_1 = \C\{\,x\mu(z) + x \otimes z\mid x\in M,\,\, z\in Z(A)_2\}
\,\subset\, T_R(M)_3 \oplus (M\otimes Z(A)_2)$$ and
$$O_2 = \C\{\,\mu(z)x + x \otimes z\mid x\in M,\,\, z\in Z(A)_2\}
\,\subset\, T_R(M)_3 \oplus (M\otimes Z(A)_2).$$
Then the kernel of the quotient map $\wt{T_R(M^*)}_3 = T_R(M^*)_3\oplus (M^*\otimes U^!) \to \tilde A^!_3$
is equal to
$$(Q_2 + O_2)^\perp + (Q_1 + O_1)^\perp= \Big((Q_2 + O_2)\cap (Q_1 + O_1)\Big)^\perp.$$

Consider any pair of elements $x\in M$ and $z\in Z(A)_2$.
Since $\mu(z)$ reduces to a central
element of $A$, there exist elements $x_1,\ldots, x_k, y_1,\ldots, y_\ell\in M$
and $r_1,\ldots, r_k, s_1,\ldots,s_\ell\in W \subset M \otimes M$
such that 
\begin{equation}\label{r-and-s}
[\mu(z),x]=\sum_{i=1}^k x_ir_i + \sum_{j=1}^\ell s_j y_j.
\end{equation}
(To avoid unwanted cancellations, we choose our elements in a way that minimizes $k+\ell$,
and we make a choice once and for all for each pair $(x,z)\in M\otimes Z(A)_2$.)
Now consider the element 
$$\kappa(x,z):=\mu(z)x-\sum s_jy_j + x\otimes z=x\mu(z)+\sum x_ir_i + x\otimes z \in (Q_2 + O_2)\cap (Q_1 + O_1).$$
The vector space $(Q_1 + O_1)\cap (Q_2 + O_2)$ is spanned by the elements 
$\kappa(x,z)$ and the subspace $Q_1\cap Q_2$, hence
any class in $T_R(M^*)_3\oplus (M^*\otimes U^!)$
that reduces to zero in $A^!$ and is orthogonal to all of the elements $\kappa(x,z)$ 
must also reduce to zero in $\tilde A^!$.
The commutator $[\eta^!, x^!]$ clearly reduces to zero in $A^!$, since $\eta^!$ is a lift of $z^!$,
which was chosen to be central.
Thus it remains only to show that it is orthogonal to each $\kappa(x,z)$.

\begin{lemma}\label{rel-prod}
If $x\in e_\a M e_\b$ and $\a<\b$, then
$$\langle z^!,j_{\a\b}(z)\rangle \, x=\sum_{i=1}^k \langle \mu^!(z^!), r_i\rangle \, x_i\in M
\,\,\,\,\text{and}\,\,\,\, s_j\in e_\a M \ep_\a M\ep_\b \,\,\text{for all}\,\, 1\leq j \leq \ell.$$
If $x\in e_\a M e_\b$ and $\b<\a$, then
$$\langle z^!,j_{\a\b}(z)\rangle \, x=\sum_{j=1}^\ell \langle \mu^!(z^!), s_j\rangle \, y_j\in M
\,\,\,\,\text{and}\,\,\,\, r_i\in \ep_\a M \ep_\b M e_\b \,\,\text{for all}\,\, 1 \leq i \leq k.$$
\end{lemma}

\begin{proof}
We will prove only the case where $\a<\b$ (the proof of the opposite case is identical).
From the definition of $\mu(z)$ in Equation \eqref{munu}, we have
$$[\mu(z), x] \in e_\a M \ep_\a M e_\a M e_\b + e_\a M e_\b M \ep_\b M e_\b
\subset e_\a M \ep_\a M (e_\a + \ep_\b) M e_\b.$$
It follows that $s_j \in e_\a M \ep_\a M (e_\a + \ep_\b)$ for all $j$.
For any $j$, the element $s_j e_\a$ lies in the subspace $\iota(W)\,\cap\, e_\a M \ep_\a M e_\a$, which
is trivial by Lemma \ref{nu}.  Hence we have $s_j \in e_\a M \ep_\a M \ep_\b$ as claimed.

Consider the specialization $\tilde A_{z^!}$ of $\tilde A$.  Since $A$ is quadratic, the natural
map from $M$ to $\tilde A_{z^!}$ is an inclusion, thus we may
regard $M$ as a subspace of $\tilde A_{z^!}$.
The image of $[\mu(z),x]$ in $\tilde A_{z^!}$ lies in $M$, and
Equation \eqref{r-and-s} tells us that it is equal to the element
$$\sum_{i=1}^k \langle z^!,\, r_i\rangle\, x_i + \sum_{j=1}^\ell \langle z^!,\, s_j\rangle\, y_j
= \sum_{i=1}^k \langle z^!,\, r_i\rangle\, x_i,$$
where the second equality follows from the fact that $e_\ga s_j e_\ga = 0$ for all $j$ and all $\ga\in\cI$.
On the other hand, it is also equal to $\langle z^!, j_{\a\b}(z)\rangle\, x$
by Proposition \ref{jabprop}, thus we obtain the desired identity.
\end{proof}

We now use Lemma \ref{rel-prod} to show that $\Big\langle [\eta^!, x^!],\, \kappa(x,z)\Big\rangle = 0$
for all $x\in e_\a M e_\b$, $x^!\in e_\b M^* e_\a$, and $z\in Z(A)_2$, and thus complete the proof of Theorem \ref{bang-flexible}.  
By Lemma \ref{noloops}, we may assume that either $\a<\b$ or $\b<\a$.
We have
\begin{align*}
\Big\langle [\eta^!, x^!],\, \kappa(x,z)\Big\rangle
&=\big\langle\mu^!(z^!)x^!,\kappa(x,z)\big\rangle - \big\langle x^!\mu^!(z^!),\kappa(x,z)\big\rangle
+ \big\langle x^!\otimes j_{\b\delta}^*(z^!)-x^!\otimes j_{\a\delta}^*(z^!),\, x\otimes z\big\rangle\\
    &= \big\langle\mu^!(z^!)x^!, x\mu(z)+\sum x_ir_i\big\rangle - \big\langle x^!\mu^!(z^!),\mu(z)x -\sum s_iy_i\big\rangle + \big\langle x^!, x\big\rangle\cdot
    \big\langle j_{\b\a}^*(z^!),\, z\big\rangle\\
    &=\big\langle\mu^!(z^!)x^!, \sum x_ir_i\big\rangle+ \big\langle x^!\mu^!(z^!), \sum s_jy_j\big\rangle + \big\langle x^!, x\big\rangle\cdot
    \big\langle z^!,\, j_{\b\a}(z)\big\rangle\\
    &=\sum \big\langle x^!, x_i\big\rangle\cdot\big\langle\mu^!(z^!), r_i\big\rangle
    + \sum \big\langle x^!, y_i\big\rangle\cdot\big\langle\mu^!(z^!), s_i\big\rangle
    - \big\langle x^!, x\big\rangle\cdot \big\langle z^!,\, j_{\a\b}(z)\big\rangle.
  \end{align*}

First assume that $\a<\b$.  By Lemma \ref{rel-prod},
each $s_i$ pairs to zero with any loops, thus the second term of the last line vanishes.
The lemma also tells us that the first and third terms cancel, 
so the entire expression is equal to zero.
Similarly, if $\b<\a$, the first term vanishes and the second and third terms cancel.
\end{proof}

Suppose that $A$ and $A^!$ are flexible.  For all $\z\in Z(\tilde A)_2$ and $\z^!\in Z(\tilde A^!)_2$,
let $\pi(\z)$ and $\pi^!(\z^!)$ denote their images in $Z(A)_2$ and $Z(A^!)_2$, respectively.
For all $\a\in\cI$, the splitting $h_\a$ of the exact sequence \eqref{exact}
and the analogous splitting $h_\a^!$ on the dual side induce a perfect pairing
$$\langle\, ,\rangle_\a:Z(\tilde A)_2\times Z(\tilde A^!)_2\to\C$$
given by the formula
$$\langle \z, \z^!\rangle_\a := \langle h_\a(\z) , \pi^!(\z^!)\rangle + \langle \pi(\z) , h^!_\a(\z^!) \rangle,$$
where we once again exploit the fact that $h_\a(z)\in U\cong Z(A^!)_2^*$ and
$h^!_\a(\z^!) \in U^!\cong Z(A)_2^*$.

\begin{proposition}\label{gmdu}
All of these pairings coincide.
\end{proposition}

\begin{proof} By definition of $j_{\a\b}$, we have 
\begin{eqnarray*}
  \langle \z, \z^!\rangle_\a-\langle \z, \z^!\rangle_\b &=&
  \big\langle j_{\b\a}\circ\pi(\z) ,\,\, \pi^!(\z^!)\big\rangle 
  + \big\langle\pi(\z) ,\,\, j^!_{\b\a}\circ\pi^!(\z^!) \big\rangle\\
&=& \big\langle \pi(\z) ,\,\, j_{\b\a}^*\circ\pi^!(\z^!)\big\rangle 
- \big\langle\pi(\z) ,\,\, j^!_{\a\b}\circ\pi^!(\z^!) \big\rangle,
\end{eqnarray*}
which vanishes by Theorem \ref{bang-flexible}.
\end{proof}

\begin{example}\label{aoo}
We illustrate Proposition \ref{gmdu}
for the algebra $A = A_{11}$ from Example \ref{two nodes}.
Though the dual algebra $A^!$ is isomorphic to $A$,
we will use separate notation in order to keep track of the two sides.
The algebra $A$ is generated by $x\in e_1 A e_2$ and $y\in e_2 A e_1$,
which satisfy the relation $yx = 0$.  Its deformation $\tilde A$ is generated
by $x, y$, and a central variable $u$, which satisfy the relation $yx = ue_2$.
On the dual side, $A^!$ is generated by $x^!\in e_2 A^! e_1$ and $y^!\in e_1 A^! e_2$,
subject to the relation $y^!x^! = 0$.  Its deformation $\tilde A^!$ is generated by $x^!,y^!$,
and a central variable $u^!$, which satisfy the relation $y^!x^! = u^!e_1$.
The generator $u\in U$ pairs to 1 with the generator $x^!y^!\in Z(A^!)_2$,
while $u^!\in U^!$ pairs to 1 with $xy\in Z(A)_2$.

The vector space $Z(\tilde A)_2$
is spanned by the elements $xy + ue_2$ and $ue_1 - xy$, and we have
$$(h_1\oplus h_2)(xy + ue_2) = (0,u)\,\,\,\,\,\text{and}\,\,\,\,\, (h_1\oplus h_2)(ue_1 - xy) = (u,0).$$ 
On the dual side, $Z(\tilde A^!)_2$ is spanned by $x^!y^! + u^!e_1$ and $u^!e_2 - x^!y^!$,
and we have $$(h_1^!\oplus h_2^!)(x^!y^! + u^!e_1) = (u^!, 0)\,\,\,\,\,\text{and}\,\,\,\,\,
(h_1^!\oplus h_2^!)(u^!e_2 - x^!y^!) = (0, u^!).$$
Let $\z = a(xy + ue_2) + b(ue_1 - xy)\in Z(\tilde A)_2$
and $\z^! = a^!(x^!y^! + u^!e_1) + b^!(u^!e_2 - x^!y^!) \in Z(\tilde A^!)_2$.
Then $$\langle \z, \z^!\rangle_1
= \langle bu, (a^!-b^!) x^!y^!\rangle +  \langle a^!u^!, (a-b)xy\rangle
= b(a^! - b^!) + a^!(a - b)
= aa^! - bb^!,$$
and $$\langle \z, \z^!\rangle_2
= \langle au, (a^!-b^!) x^!y^!\rangle + \langle b^!u^!, (a-b)xy\rangle
= a(a^! - b^!) + b^!(a - b)
= aa^! - bb^!.$$
Thus the two pairings are the same.
\end{example}

We may now use Proposition \ref{gmdu} to prove Theorem \ref{second}, which we restate here.

\begin{corollary}\label{main}
If $A$ is flexible, then 
$\cZ(A)$ and $\cZ(A^!)$ are canonically GM dual.
\end{corollary}

\begin{proof}
By Definition \ref{fibered}, we associate to the localization algebra $\cZ(A)$ the fibered
arrangement consisting of the subspaces $$H_\a := h_\a^*(U^*)\subset Z(\tilde A)^*_2,$$
each of which projects isomorphically onto $U^*$.
Definition \ref{duality} tells us that a duality between $\cZ(A)$ and $\cZ(A^!)$
is a perfect pairing between $Z(\tilde A)^*_2$ and $Z(\tilde A^!)^*_2$ such that
the kernels of the two projections are perpendicular to each other, as are $H_\a$
and $H_\a^!$ for each $\a\in\cI$.

For each $\a\in\cI$, we have constructed a perfect pairing
$$\langle\, ,\rangle_\a:Z(\tilde A)_2\times Z(\tilde A^!)_2\to\C,$$
which induces a dual pairing 
$$\langle\, ,\rangle^*_\a: Z(\tilde A)^*_2\times Z(\tilde A^!)^*_2\to\C.$$
It is clear from the definition of the pairing
that the kernels of the two projections are perpendicular spaces of each other,
and that $H_\a$ 
is the perpendicular space to $H_\a^!$.
By Proposition \ref{gmdu}, the pairings $\langle\, ,\rangle_\a$ all coincide,
therefore we have one canonical pairing satisfying all of the required properties.
\end{proof}

\section{Example: Polarized arrangements and hypertoric varieties}\label{hypertoric examples}
In this section and the next, we consider two families of examples of Koszul 
dual pairs of flexible algebras,
along with the associated dual pairs of localization algebras.  As we will see, most of our examples
have cohomological interpretations in addition to algebraic ones.
For our first example, we use a ring that we introduced in 
an earlier paper \cite{GDKD}, constructed from the following linear algebra data.

\begin{definition}\label{polarized}
A {\bf polarized arrangement} $\cV$ is a triple $(V, \eta, \xi)$, where
$V$ is a linear subspace of a coordinate vector space $\R^n$,
$\eta\in\R^n/V$, and $\xi\in V^*$.
\end{definition}

It is convenient to think of these data as describing
an affine space $V_\eta \subs \R^n$ given by translating 
$V$ away from the origin by $\eta$, together with an affine linear
functional on $V_\eta$ given by $\xi$ and a finite 
hyperplane arrangement $\cH$ in $V_\eta$, whose
hyperplanes are the (possibly empty) restrictions of the coordinate hyperplanes in $\R^n$.
We will assume
that $\eta$ and $\xi$ are chosen generically enough so that
$\cH$ is simple (any set of $m$ hyperplanes intersects in codimension $m$
or not at all) and $\xi$ is non-constant on any positive-dimensional intersection
of the hyperplanes.
For now (until Remark \ref{nongeometric}) we will also assume that $\cV$ is {\bf rational}, meaning that
$V$, $\eta$, and $\xi$ are all defined over $\Q$.

In \cite[\S 4]{GDKD}, we explained how to associate to this data a standard Koszul algebra $B(\cV)$.
We sketch this construction here; many more details are given in \cite{GDKD}.
For all $\a\in\{\pm 1\}^n$, let 
$$\Delta_\a = \big\{v\in V_\eta\subset\R^n\,\,\big{|}\,\, \text{$\a(i)\cdot v_i\geq 0$ for all $i = 1,\ldots n$}\big\}.$$
Geometrically, $\Delta_\a$ is the chamber of $\cH$ consisting of vectors that lie on a fixed
side of each hyperplane.
Let the indexing set $\cI$ be the set of sign vectors $\a$ such that $\Delta_\a$ is nonempty and 
the affine linear functional $\xi$ is bounded above on $\Delta_\a$.
To each $\a\in\cI$, we may associate a toric variety $X_\a$, with an effective action of the algebraic
torus $T$ whose Lie algebra is equal to $V^*_\C$ and whose character lattice is
$\Z^n \cap V \subset V_\C$.  The action of the maximal compact subtorus
is hamiltonian, and $\Delta_\a$ is the moment polyhedron for this action.

For all $\a,\b\in\cI$, let $d_{\a\b}$ be 
codimension of $\Delta_\a\cap\Delta_\b$ in $V_\eta$, and let $X_{\a\b}$ be the toric variety
with moment polyhedron $\Delta_\a\cap\Delta_\b$.
As a graded vector space, $B(\cV)$ is defined as the sum
\begin{equation}\label{bdef}
\bigoplus_{\a,\b\,\in\,\cI} H^*(X_{\a\b})[-d_{\a\b}].
\end{equation}
The product that we define is a convolution product:  to multiply an element of
$H^*(X_{\a\b})$ with an element of $H^*(X_{\b\ga})$, we pull both classes
back to the toric variety with moment polyhedron $\Delta_\a\cap\Delta_\b\cap\Delta_\ga$,
multiply them there, and then push forward to $X_{\a\ga}$.  It is an easy combinatorial exercise
to check that this product respects the grading.  Showing that it is associative is more subtle,
and in fact is only true if we push forward not with respect to the complex orientations,
but with respect to a collection of combinatorially defined orientations on the various toric varieties \cite[4.10]{GDKD}.

\begin{remark}\label{hyp-conj}
The motivation for this definition comes from the geometry of 
the {\bf hypertoric variety} $\MV$ associated to $\cV$, which 
is a complex symplectic algebraic variety of 
dimension $2\dim V$ (or a hyperk\"ahler manifold of real dimension $4\dim V$).
It comes equipped with an effective hamiltonian action of $T$ 
(or a tri-hamiltonian action of the maximal compact subtorus).
The variety itself depends only on $\cH$, and the covector $\xi$ determines an action of $\C^*$
on $\MV$.  For each $\a\in\cI$, the toric variety $X_\a$
sits inside of $\MV$ as a Lagrangian subvariety.  The union of all of these subvarieties
is equal to the set of points $p\in\MV$ such that $\lim_{\la\to\infty}\la\cdot p$ exists.
We conjecture that the algebra $B(\cV)$ is isomorphic to the Ext-algebra 
in the Fukaya category of $\MV$
of the sum of the objects associated to the Lagrangian subvarieties $X_\a$.
For more information on hypertoric varieties, see the survey article \cite{Pr07}.
\end{remark}

Given a polarized arrangement $\cV$, we define its {\bf Gale dual}
$\cV^\vee = (V^\perp, -\xi, -\eta)$, where $V^\perp$ sits inside of the dual coordinate
vector space $(\R^n)^*$, 
$-\xi\in V^*\cong (\R^n)^*/V^\perp$, and $-\eta \in \R^n/V\cong (V^\perp)^*$.

\begin{theorem}\label{bv}{\em\cite[3.11, 4.14, 4.16, 5.23, \& 5.24]{GDKD}}
The algebra $B(\cV)$ is standard Koszul, and its center is isomorphic as a graded ring to the cohomology ring of $\MV$.  The algebras $B(\cV)$ and $B(\cV^\vee)$ 
are Koszul dual to each other.
\end{theorem}

\begin{remark}
Note that for $B(\cV)$ and $B(\cV^\vee)$ to be Koszul dual, their
degree 0 parts must be isomorphic.  The degree 0 part of $B(\cV)$
is spanned by the unit elements $1_{\a\a} \in H^*(X_{\a\a})$ for all $\a\in\cI$.
Let $\cI^\vee$ be the corresponding set for $\cV^\vee$, that is, the set of sign
vectors that give chambers of $\cV^\vee$ on which $-\eta$ is bounded above.
We prove in \cite[2.4]{GDKD} that $\cI^\vee = \cI$, and therefore that
there is a canonical isomorphism between $B(\cV)_0$ and $B(\cV^\vee)_0$.
\end{remark}

In \cite{GDKD} we also define a deformation $\tilde{B}(\cV)$ of
$B(\cV)$. 
It is defined by replacing all of cohomology rings
in (\ref{bdef}) with $T$-equivariant cohomology rings:
\begin{equation}\label{btdef}
\tilde{B}(\cV)\,\,:= \bigoplus_{\a,\b\,\in\,\cI} H^*_T(X_{\a\b})[-d_{\a\b}],
\end{equation}
with a convolution product defined as for $B(\cV)$.
By $\cite[4.5 \;\&\; 4.10]{GDKD}$
it is a flat deformation of $B(\cV)$ over $V^*_\C$, where 
$\tilde{B}(\cV) \to B(\cV)$ is the map forgetting the equivariant
structure, and the map
\[S := \Sym(V_\C) \cong H^*_T(pt) \to \tilde{B}(\cV)\] 
sends an element of $S$ to the sum of its
images in $H^*_T(X_{\a\a})$ over all $\a \in \cI$.

\begin{proposition}\label{eqb}  The deformation $\tilde{B}(\cV)$ is 
flexible and malleable. Its center
(with localization algebra structure defined in Corollary \ref{gma})
is isomorphic as a localization algebra to the $T$-equivariant cohomology ring of $\MV$
(with localization algebra structure defined in Example \ref{geom-gm}).
\end{proposition}

\begin{proof}
The isomorphism of $S$-algebras 
between the center of $\tilde B(\cV)$ and $H^*_T(\MV)$
is given in \cite[4.16]{GDKD}, where 
we show that both rings are quotients
of the polynomial ring $\C[u_1,\ldots,u_n]$ by the same ideal.
This result also shows that $Z(\tilde{B}(\cV)) \to Z(B(\cV))$ is
surjective, so $\tilde{B}(\cV)$ is a flexible deformation.
We also exhibit in \cite[\S 2.6]{GDKD} a natural
bijection between $\cI$ and the fixed point set $\MV^T$; it sends
$\a$ to the fixed point $x_\a \in X_{\a}^T \subset \MV^T$ corresponding
to the vertex of $\Delta_\a$ on which $\xi$ attains its maximum.

The standard modules over $B(\cV)$ are described
geometrically by \cite[5.22]{GDKD}; we have
\[V_\a \cong \bigoplus_{\b \in \cI} H^*(\{x_\a\} \cap X_{\b})[-d_{\a\b}],\]
with a natural right action of $B(\cV)$ by convolution.
Corollary \ref{deformed standard characterization} now implies that
\[\tilde{V}_\a := \bigoplus_{\b \in \cI} H^*_T(\{x_\a\} \cap X_{\b})[-d_{\a\b}],\]
with the action of $\tilde{B}(\cV)$ by convolution, is the deformed standard object defined 
in Section \ref{sec-malleable}.  It follows immediately
that the map $h_\a\colon Z(\tilde{B}(\cV))_2 \to V_\C$ 
of Proposition \ref{hadef} coincides
with the localization map $H^2_T(\MV) \to H^2_T(x_\a)$.

To see that $\tilde{B}(\cV)$ is malleable, first note that the anti-involution
induced by the isomorphism of $X_{\a\b}$ with $X_{\b\a}$
induces an isomorphism $B(\cV)\cong B(\cV)^\text{op}$.
Malleability now follows from the fact that
the localization map $H^*_T(\MV)\to H^*_T\!\left(\MV^T\right)$ is an isomorphism
over the generic point of $\mathfrak{t}$.
\end{proof}

Using this, we have a simple description of 
the degree two part of the maps $h_\a$ from the localization algebra structure, and hence
of the associated fibered arrangement.  For simplicity we will
assume that $V$ is not contained in any coordinate hyperplane, 
so there are no empty hyperplanes in our arrangement. 

For each $\a \in \cI$, let $p_\a\in\Delta_\a$ be the point at which $\xi$ attains its maximum and
let $b_\a \subset \{1,\dots, n\}$ be the set of indices $i$ for which 
the $i^\text{th}$ hyperplane of $\cH$ contains $p_\a$.
The collection $\{b_\a\mid\a\in\cI\}$ consists precisely of all subsets of $\{1,\dots, n\}$
for which then the composition of the inclusion $\iota\colon V \hookrightarrow \R^n$
with the coordinate projection $\pi_\a\colon \R^n \to \R^{b_\a}$
is an isomorphism.  Such subsets are known as the {\bf bases} of $\cV$.

\begin{proposition}\label{hypertoric localization algebra}
There is an isomorphism of $H^2_T(\MV)$ with $\C^n$ such that the
inclusion $$V_\C\cong H^2_T(pt) \hookrightarrow H^2_T(\MV)\cong\C^n$$ is the 
complexification of $\iota$.  Under this identification, the
restriction of the localization algebra map $h_\a$ to degree $2$ 
is the complexification of $(\pi_\a \circ \iota)^{-1} \circ \pi_\a$.
The fibered arrangement associated to $\tilde{B}(\cV)$ is the union
of the coordinate subspaces $(\C^{b_\a})^*$ of the dual space 
$(\C^n)^* \cong H_2^T(\MV)$.
\end{proposition}

\begin{proof}  The first statement follows from the standard 
description of the equivariant cohomology of a hypertoric
variety \cite[3.2.2]{Pr07}.  The remaining statements follow
easily from the fact that $\C^{b_\a^c}$ must be in 
the kernel of $h_\a$.
\end{proof}

\begin{theorem}\label{hypertoric universal} Suppose that the subspace
$V\subset \R^n$ contains no coordinate line.
Then $\tilde{B}(\cV)$ is isomorphic to the universal
deformation of $B(\cV)$. 
\end{theorem} 
\begin{proof}
Since $V$ contains no coordinate line, 
$V^\bot$ is not contained in any coordinate
plane, and so Theorem \ref{bv} implies that 
$$Z(B(\cV)^!)_2 \cong Z(B(\cV^\vee))_2 \cong H^2(\mathfrak{M}(\cV^\vee)) \cong V_{\C}^*,$$
where the last isomorphism comes from the formula for the 
cohomology of a hypertoric variety \cite{Ko, HS, Pr07}, which
gives $H^2(\mathfrak{M}(\cV^\vee)) \cong \C^n/V_{\C}^\bot = V_{\C}^*$.

Thus to prove that the map $\psi\colon V_\C \to Z(B(\cV)^!)_2^*$
associated to the deformation $\tilde{B}(\cV)$ by Theorem \ref{flat}
is an isomorphism, it is enough to show that it is surjective.
We can do this using Lemma \ref{image of psi}.  Using  
Proposition \ref{hypertoric localization algebra} and a little
linear algebra, it is not hard to show that for 
any $\a, \b \in \cI$, the value of $j_{\a\b}$ on
the parameter $\eta \in \C^n/V_\C \cong Z(B(\cV))_2$
is given by
\[j_{\a\b}(\eta) = p_\b - p_\a \in V\subset V_\C \cong Z(B(\cV)^!)_2^*.\]
(Note that $p_\a$ and $p_\b$ both lie in the affine space $V_\eta$, 
so their difference lies in the vector space $V$.)
The surjectivity of $\psi$ now follows from 
Lemma \ref{image of psi}, using the fact
that the points $p_\a$ form an affine spanning set for $V_\eta$ (this is where we use the
assumption that $V$ contains no coordinate line).
\end{proof}

\begin{corollary}\label{htgm} The 
localization algebras $\cZ(\tilde{B}(\cV))$ and $\cZ(\tilde{B}(\cV^\vee))$
are dual.  
\end{corollary}

Corollary \ref{htgm} follows immediately from Corollary \ref{main}, but
the concrete description of these 
fibered arrangements in Proposition \ref{hypertoric localization algebra}
makes it easy to see this directly: if 
$b^\vee_\a \subset \{1, \dots, n\}$ is the basis for $\cV^\vee$ indexed by $\a$, 
then we have $b^\vee_\a = b_\a^c$ \cite[2.9]{GDKD}, and so
$H^\vee_\a=\C^{b^\vee_\a}$ is perpendicular to $(\C^{b_\a})^*=H_\a$.

\begin{example}\label{firstgen}
  If $V$ is $(n-1)$-dimensional and $\cH$ consists of a collection of
  $n$ hyperplanes in general position, then the hypertoric variety
  $\MV$ is isomorphic to the cotangent bundle of $\mathbb{P}^{n-1}$.
  Dually, if $V$ is one-dimensional and $\cH$ consists of $n$ points
  on a line, then $\MV$ is isomorphic to the minimal resolution of the
  symplectic surface singularity $\C^2/\Z_n$, which retracts onto a
  chain of $n-1$ projective lines.  Thus the duality of 
  $\cZ(\tilde{B}(\cV))$ and $\cZ(\tilde{B}(\cV^\vee))$ 
  generalizes that of Examples \ref{ptwo} and \ref{two lines}.
\end{example}

\begin{remark}\label{nongeometric}
  We have assumed that $\cV$ is rational in order to give that shortest
  and best motivated definition of the algebra $B(\cV)$.  In
  \cite{GDKD}, however, we do not make this assumption.  
  Although the toric varieties $X_\a$ and the
  hypertoric variety $\MV$ are not defined when $\cV$ is not rational,
  it is still possible to give combinatorial definitions of rings and
  localization algebras that specialize to the ordinary and equivariant
  cohomology rings of these spaces in the rational case.  In this more
  general setting, Theorem \ref{bv}, Proposition \ref{eqb}, Proposition
  \ref{hypertoric localization algebra}, Theorem \ref{hypertoric universal}, and Corollary \ref{htgm}
  go through
  exactly as stated.  The rings $B(\cV)$ for $\cV$ not rational are
  the only examples that we know of flexible algebras that are not
  associated with any algebraic variety.
\end{remark}

\section{Example: Category $\mathcal{O}$ and Spaltenstein varieties}\label{category O examples}
In this section, we apply our deformation result to integral 
blocks of parabolic 
category $\cO$ for the Lie algebra $\mg = \mathfrak{gl}_n(\C)$.
In particular, we show that the universal deformation of the 
endomorphism algebra of a minimal
projective generator of such a block is malleable and flexible,
and identify the associated localization algebra, which turns out to come from the
equivariant cohomology ring of a Spaltenstein variety.  
We accomplish this by identifying modules over our universal
deformation with objects in ``deformed category $\cO$'' as considered
by Soergel and Fiebig.
Most of the results of this section should be true for $\mg$ an arbitrary reductive 
Lie algebra, but in our proofs we use Brundan's computation of centers of blocks of 
parabolic category $\cO$ for $\mathfrak{gl}_n$ and therefore restrict ourselves to $\mg=\mathfrak{gl}_n(\C)$.

We mostly follow Brundan's notation to describe the
blocks of $\cO$.  Let $\mg = \mathfrak{gl}_n(\C)$, and let 
$\mb$ and $\mh$ be the Lie subalgebras
of $\mg$ consisting of upper triangular and diagonal
matrices, respectively.
The Weyl group $W = W(\mg, \mh)$ is the symmetric group 
$S_n$, which acts on $\mh^* \cong \C^n$ and on the 
weight lattice $\Lambda = X(T) \subset \mh^*$ by permuting coordinates.
With these conventions, a weight $\lambda\in \Lambda$ 
is dominant if and only if $\lambda_i > \lambda_j$ for all $i < j$.
Let $w_0 \in W$ denote the longest element.

A {\bf composition} of $n$ is a doubly-infinite
sequence $\nu = (\dots, \nu_{-1}, \nu_0, \nu_1, \nu_2, \nu_3, \dots)$ of non-negative
integers whose sum is $n$.  
Given such a $\nu$, there is a unique dominant weight $\a_\nu$ with $\nu_i$
entries equal to $-i$ for every integer $i\in\Z$.
This gives a bijection between the $W$-orbits $\Lambda/W$ 
of the weight lattice and the set of compositions of $n$.

For any composition $\nu$ of $n$, let
$W_\nu = W_{\alpha_\nu} \subset W$ be the 
stabilizer of the dominant weight $\alpha_\nu \in \Lambda$.  
Elements of $W_\nu$
are permutations which preserve subsets of consecutive elements of 
$\{1, \dots, n\}$ of sizes $\dots \nu_{-1}, \nu_0, \nu_1, \dots$, 
in that order.  We will refer to these subsets as \emph{$\nu$-blocks}.
We also define another associated composition $\bar{\nu}$ by letting 
$\bar\nu_1 = \nu_i$ where $i$ 
is the smallest index with $\nu_i \ne 0$, $\bar\nu_2 = \nu_{i'}$
where $i'$ is the next smallest index where $\nu$ is nonzero, 
and so on, letting all other $\bar\nu_j$ be zero.  

Associated to $\mg = \mathfrak{gl}_n$ we have the
Bernstein-Gelfand-Gelfand category $\cO$ of all finitely 
generated $\fg$-modules which are $\mh$-diagonalizable
and locally finite over $\mb$.  For simplicity we add
the additional assumption that
all weights lie in the lattice $\Lambda$.  For $\alpha\in \Lambda$,
there is a unique simple $\mg$-module $L(\a)$ with highest weight
$\a - \rho$, where $\rho = (0, -1, -2, \dots, -n + 1)$.  
These are the simple objects of $\cO$.

For a composition 
$\nu$ of $n$, define $\cO_\nu$ to be
the Serre subcategory of $\cO$ generated by all $L(w\alpha_\nu)$
for $w \in W$.  Obviously the weight
$w\alpha_\nu$ only depends on the image of $w$ in 
$W/W_\nu$. Define a composition $\nu^t$ 
by letting $\nu^t_j$ be the number of $i\in\Z$
for which $\nu_i \ge j$ if $j\geq 1$, and zero otherwise.  It is a {\bf partition}, meaning 
a composition that's supported on $\mathbb{N}$ and non-increasing.  The partition
$\nu^+ := (\nu^t)^t$ has the same parts as 
$\nu$, sorted into non-increasing order.

Given another composition $\mu$ of $n$, let
$\mp = \mp_\mu$ be the parabolic subalgebra
of $\mg$ given by all block upper-triangular 
matrices where the blocks are the $\mu$-blocks.
Associated to $\mu$ we have the parabolic category
$\cO^\mu = \cO^\mp$, the full subcategory
of $\cO$ of objects which are $\mp$-locally finite.
Its simple objects are $\{L(\alpha)\mid\a\in\Lplusmu\}$, where
\[\Lplusmu := 
\{\alpha \in \Lambda \mid \alpha_j > \alpha_k\;\text{whenever $j<k$ and $j, k$ lie in the same $W_\mu$-orbit}\}.
\] 
Define $\cO^\mu_\nu := \cO^\mu \cap \cO_\nu$;
it is the Serre subcategory of $\cO$ generated 
by the simple objects $L(w\alpha_\nu)$ for $w \in \cI^\mu_\nu$, 
where 
\[\cI^\mu_\nu = \{w \in W/W_\nu \mid w\alpha_\nu \in \Lplusmu\}.\] 

\begin{lemma}
The category $\cO^\mu_\nu$ is nonzero if and only if
$\mu^+ \le \nu^t$ in the dominance order
on partitions.  The map $w \mapsto W_\mu w$
defines a bijection between $\cI^\mu_\nu$
and the set of double cosets in 
$W_\mu \backslash W/W_\nu$ of size
$|W_\mu| \times |W_\nu|$.
\end{lemma}

The category $\cO^\mu_\nu$ has enough projectives.
For $w \in \cI^\mu_\nu$, let $P^\mu(w\a_\nu)$ be a
projective cover of $L(w\a_\nu)$ in $\cO^\mu_\nu$.
(Note that although $L(w\a_\nu)$ can lie in $\cO^\mu_\nu$
for many choices of $\mu$, in general the projective
covers in these categories will be different.)
Let $P^\mu_\nu := \bigoplus_{w \in \cI^\mu_\nu} P^\mu(w\a_\nu)$
be a minimal projective generator of $\cO^\mu_\nu$, 
and let $A^\mu_\nu := \End(P^\mu_\nu)$, so that 
$M \mapsto \Hom(P^\mu_\nu, M)$ is an equivalence
between $\cO^\mu_\nu$ and the category of finitely
generated right $A^\mu_\nu$-modules.

For a composition $\mu$, let $\mu^o$ 
denote the reversed composition given by $\mu^o_{i} := \mu_{- i}$.

\begin{proposition}\label{KD for cat O}
 The ring $A^\mu_\nu$ has a grading with respect to which it 
is standard Koszul.  There is an 
isomorphism
$(A^\mu_\nu)^! \cong A^{\nu}_{\mu^o}$
whose map on idempotents is induced by the
map $\cI^\mu_\nu \to \cI^\nu_{\mu^o}$ taking 
$wW_\nu$ to $w^{-1}w_0 W_{\mu^o}$, where
$w_0$ has maximal length in the coset $wW_\nu$.
\end{proposition}
\begin{proof}
The construction of a Koszul grading and the identification of
the Koszul dual is accomplished in Backelin \cite[1.1]{Back99}.  
The fact that $A^\mu_\nu$ is quasi-hereditary follows from
\cite[Theorem 6.1]{RC}.  
Since both $A^\mu_\nu$ and its dual are
quasi-hereditary, and the associated partial orders on 
the idempotents are reversed, 
\cite[Theorem 3]{ADL03} implies that $A^\mu_\nu$ is standard Koszul. 
\end{proof}

\begin{remark}
Up to isomorphism the algebra $A^\mu_\nu$
only depends on the subgroups $W_\nu, W_\mu \subset W$.
This is obvious for $\mu$, while for $\nu$
the required equivalences are given by 
translation functors.  In addition, the rings
$A^\mu_\nu$ and $A^{\mu^o}_{\nu^o}$ are
isomorphic, using the automorphism of $\mg$
given by the adjoint action by any representative
for $w_0$ in $G = GL_n(\C)$.
\end{remark}

Since $A^\mu_\nu$ is Koszul, we can consider our
universal graded deformation $\tilde{A}^\mu_\nu$
and its associated localization algebra $\cZ(\tilde{A}^\mu_\nu)$
as given by Corollary \ref{gma}.  We wish to 
relate this localization algebra with one arising from
geometry, specifically, from the equivariant
cohomology of a Spaltenstein variety.

As before, let $\mu$, $\nu$ be compositions of $n$,
and suppose that $\mu^+ \le \nu^t$.  Let
$P_\nu \subset G := \operatorname{GL}_n(\C)$ be the parabolic subgroup
with Lie algebra $\mp_\nu$, and let
\[X_\nu = G/P_\nu = \{0 = F_0 \subset F_1 \subset F_2 \subset \dots 
\subset \C^n \mid \dim_\C F_i = \bar{\nu}_1 + \dots + \bar{\nu}_i\}
\] be the associated partial 
flag variety.
(Note that on the geometric side it is $\nu$, not $\mu$,
which specifies the parabolic; this is related to Remark \ref{langlands dual} below.)  
The cotangent bundle $T^*X_\nu$
may be identified with the variety of pairs
$$\{(F_\bullet, N)\in X_\nu\times\mg\mid N F_i\subset F_{i-1}\;\text{for all}\; i > 0\}.$$
The moment map $\pi\colon T^*X_\nu \to \mg^*$ is a resolution of 
the closure of a nilpotent coadjoint orbit; the
Spaltenstein variety $X^\mu_\nu$ is the fiber of
this map over a point in an orbit of type $\mu$.
More precisely, we identify $\mg^*$ with $\mg$
via the inner product $\langle A, B\rangle := \operatorname{tr}(AB)$, and
we let $N_\mu \in \mg^*$ be the nilpotent matrix
defined by $N_\mu(e_i) = e_{i+1}$ if $i$ and $i+1$ are in the
same $\mu$-block, and $N_\mu(e_i)=0$ otherwise, where $e_i$ is
the $i^\text{th}$ standard basis element of $\C^n$.  Then we define
\[X^\mu_\nu := \pi^{-1}(N_\mu) = \{F_\bullet \in X_\nu \mid N_\mu F_i \subset F_{i-1} \;\text{for all}\; i > 0\}.\]
Let $T\subset G$ be the maximal torus consisting of diagonal matrices.
It acts naturally on the flag variety $X_\nu$,  and
the subtorus $T^\mu := T^{W_\mu} = Z_G(N_\mu) \cap T$
preserves the subvariety $X^\mu_\nu$.

\begin{remark}\label{langlands dual}
The group $G$ whose flag variety we have introduced
is morally the Langlands dual of the group with Lie algebra $\mg$ whose representations we are studying.
In particular, the Lie algebra $\mt$ of $T$ should be identified with the dual of $\mh$.
The Lie algebra $\mt^\mu$ is a subspace of $\mt$, so its dual $\mh^\mu$ should be thought
of as a quotient of $\mh$.

This can be confusing, since the group $\operatorname{GL}_n(\C)$ is isomorphic
to its own Langlands dual.  In particular, $\mh^\mu$ is isomorphic to $\mt^\mu$, and the quotient
map from $\mh$ to $\mh^\mu$ is given by averaging over $W_\mu$.  We will, however,
be careful never to use this isomorphism;
we will always distinguish between $\mt$ and $\mh$.
\end{remark}

Recall that the $T$-fixed points in the flag variety $X_\nu$
are in bijection with $W/W_\nu$ by $w \mapsto p_w$, where
$p_w$ is the flag $F_\bullet(w)$ given by 
$F_i(w) = \Span \{e_{w(j)} \mid 1 \le j \le \bar\nu_1 + \dots + \bar\nu_i\}$.

\begin{proposition}\label{Spaltenstein fixed points}
 The set of $T^\mu$-fixed points in $X^\mu_\nu$ is
$\{p_w \mid w \in \cI^\mu_\nu\}$.
\end{proposition}   

\begin{proof}  First we show that a $T^\mu$-fixed point
in $X^\mu_\nu$ must in fact be fixed by $T$.
Let $V_1, \dots, V_r$ be the vector spaces spanned by the 
standard basis vectors in the $\mu$-blocks, so that 
an element of $T$ lies in $T^\mu$ if and only if it
acts on each $V_j$ by multiplication by a scalar.
Let $F_\bullet \in X^\mu_\nu$ be fixed by $T^\mu$;
we need to show that each $F_i$ is spanned by 
vectors in the standard basis. 
Since $F_\bullet$ is fixed by $T^\mu$, we have 
$F_i = \bigoplus_{j = 1}^r F_i \cap V_j$ for each $i$.
The operator $N_\mu$ acts on $V_j$ as a regular nilpotent 
matrix which sends each $e_k$ to $e_{k+1}$ or to $0$.
Since $N_\mu(F_1 \cap V_j) = 0$ and $N_\mu(F_i \cap V_j) \subset F_{i-1}\cap V_j$ for $i > 1$,
it follows that $F_i \cap V_j$
is spanned by vectors $e_k$, and so $F_i$ is as well.


Thus all our fixed points are of the form $p_w$ for some $w \in W/W_\nu$.
In order to have $p_w \in X^\mu_\nu$, we must have 
$N_\mu(F_i(w)) \subset F_{i-1}(w)$ for all $i > 0$. This is equivalent
to saying that if $j < k$ lie in the same $\mu$-block, then
$w^{-1}(j)$ lies in a later $\nu$-block than $w^{-1}(k)$, or equivalently,
$(w\alpha_\nu)_j > (w\alpha_\nu)_k$.  Therefore $p_w \in X^\mu_\nu$ if and only if 
$w \in \cI^\mu_\nu$. 
\end{proof}

We will show below that the torus-equivariant cohomology of $X_\nu^\mu$ is 
isomorphic to the localization algebra of the universal deformation of $A^\mu_\nu$.  
However, the torus $T^\mu$ is too large
in general.  For instance, take $n = 3$ and let $\mu = \nu$
be the partition $(2, 1)$.  Then $\cI^\mu_\nu$ has only
one element, so $A^\mu_\nu \cong \C$ and its degree two
part is zero, while $T^\mu$ is two-dimensional.
The action of $T^\mu$ will factor through a
quotient torus which we define as follows.

Let $\lambda = \nu^t$, the transpose partition to $\nu$.
\begin{definition}
  Let $\cJ$ be the collection of all subsets $J \subset \{1, \dots,
  n\}$ such that $J$ is a union of
  $\mu$-blocks and
    \[|J| = \lambda_1 + \dots + \lambda_k,\]
    where $k$ is the number of $\mu$-blocks appearing in $J$.
 \end{definition}
 
Note that the existence of such a $J$ other than
$\{1,\ldots,n\}$ implies that there is an index where the dominance
inequality required for $\mu^+\leq \la$ is an equality.  Thus, the elements of
$J$ measure where $\mu^+$ comes closest to not being less than $\la$.  

For any $J \in \cJ$, let $\mathbf{1}_J := \sum_{i \in J} e_i \in \mt^\mu\subset\mt\cong\C^n$,
and let $T^\mu_\nu$ be the quotient of $T^\mu$ by the connected 
subtorus with Lie algebra spanned by $\{\mathbf{1}_{J}\mid J \in \mathcal{J}\}$,
so
\[\Lie T^\mu_\nu = \mt^\mu_\nu := \mt^\mu/\Span\{\mathbf{1}_{J} \mid J \in \mathcal{J}\}.\]

The meaning of the sets $J \in \cJ$ is explained by the following 
combinatorial result.  For any $J \subset \{1,\dots, n\}$, 
let $W_{\! J} := \{w \in W \mid w(J) = J\}$.
\begin{lemma} \label{extra restrictions}
 For any $J \in \cJ$, the set
$\{w\a_\nu \mid w \in \cI^\mu_\nu\}$ is contained
in a single $W_{\! J}$-orbit.  In particular, the inner
product $\langle w\a_\nu, \mathbf{1}_J\rangle$ is 
independent of $w$ for $w \in \cI^\mu_\nu$.
\end{lemma}
\begin{proof} We can assume that $\mu^+ \le \lambda$,
since otherwise $\cI^\mu_\nu = \emptyset$.  Take 
an element $J \in \cJ$, and suppose that
$|J| = \lambda_1 + \dots + \lambda_k$.  
Fix a permutation $w \in W$. 
The vector $w\alpha_\nu$ has $\nu_{i}$ entries 
equal to $-i$ for all $i \in \Z$, so if we let
$m_i = \#\{j \in J \mid (w\alpha_\nu)_j = -i\}$, we have $m_i \le \nu_i$
for all $i$.   If $w \in \cI^\mu_\nu$ then $w\alpha_\nu$
lies in $\Lplusmu$, so the entries in each $\mu$-block are strictly decreasing.
In particular, each $\mu$-block has distinct entries, and since
$J$ is the union of exactly $k$ $\mu$-blocks, we must have
$m_i \le k$ for all $i \in \Z$.
But then
\[|J| = \sum_{i \in \Z} m_i \le \sum_{i \in \Z} \min(\nu_i, k) = \lambda_1 + \dots + \lambda_k = |J|,\]
so we must have $m_i = \min(\nu_i, k)$ for all $i$.  This means that
the multiset of entries of $w\alpha_\nu$ in the places $j \in J$ is
independent of $w \in \cI^\mu_\nu$.
\end{proof}

\begin{proposition}  The action of $T^\mu$ on $X^\mu_\nu$ factors through
the quotient $T^\mu \to T^\mu_\nu$.
\end{proposition}

\begin{proof}  Take any $J \in \cJ$, and let $T_{\! J}\subset T^\mu$ be
the connected subtorus with Lie algebra $\C \cdot \mathbf{1}_J$.
We need to show that $T_{\! J}$ acts trivially on $X^\mu_\nu$.
Suppose not; then $X^\mu_\nu$ must meet more than one   
connected component of the fixed point set $(X_\nu)^{T_{\! J}}$,
and so the fixed point set 
$(X^\mu_\nu)^{T^\mu}$ must also meet more than one component.  
But two $T$-fixed points $p_w$, $p_{w'}$ lie in the
same component of $(X_\nu)^{T_{\! J}}$ if and only if 
$w$ and $w'$ lie in the same $W_{\! J}$-orbit, contradicting 
Lemma \ref{extra restrictions}.
\end{proof}

The following is our main result relating Spaltenstein varieties with 
category $\cO$.  The next section is devoted to its proof.

\begin{theorem}\label{cat O and Spaltensteins} 
The universal deformation $\tilde{A}^\mu_\nu$ of 
$A^\mu_\nu$ is malleable and flexible.  
There is an isomorphism of localization algebras
between $Z(\tilde{A}^\mu_\nu)$ and $H^*_{T_\mu^\nu}(X^\mu_\nu)$.
\end{theorem}

\begin{remark}\label{Brundan}
Theorem \ref{cat O and Spaltensteins} implies that 
the rings $Z(A^\mu_\nu)$ and $H^*(X^\mu_\nu)$ are isomorphic.  
Another proof of this can be obtained by using \cite{Bru06},
which gives a presentation by generators and relations for $Z(A^\mu_\nu)$, and \cite{BrO}, which
shows that the same presentation gives $H^*(X^\mu_\nu)$.   The
proof of Theorem \ref{cat O and Spaltensteins} uses some 
technical results of Brundan from \cite{Bru08} which 
were also used in \cite{Bru06}, so our proof is not 
entirely new.

Note also that Theorems \ref{KD for cat O} and \ref{cat O and Spaltensteins} 
imply that $\mt^\mu_\nu$ and $Z((A^\mu_\nu)^!)_2 \cong Z(A^{\nu}_{\mu^o})_2 \cong Z(A^{\nu^o}_\mu)_2$ 
must be isomorphic, as they can all be interpreted as the base of the universal
deformation of $A^\mu_\nu$.  This also follows from Brundan's work: in
fact, our description of $\mt^\mu_\nu$ is exactly the degree two 
part of the isomorphism $H^*(X^{\nu^o}_\mu) \cong Z(A^{\nu^o}_\mu)$ given in
\cite{Bru06}.
\end{remark}

\begin{remark}\label{fukaya}
Like the algebras $B(\cV)$ 
considered in Section \ref{hypertoric examples}, 
the algebras $A^\mu_\nu$ are
related to the geometry of certain symplectic algebraic varieties.
Let $S^\mu$ be the {\bf Slodowy slice} to the nilpotent matrix $N_\mu$, 
constructed explicitly in \cite[\S 7.4]{Slod80}, and recall the moment map
$\pi\colon T^*X_\nu \to \mg^*$.
The preimage $\tilde{S}^\mu_\nu := \pi^{-1}(S^\mu)$ is a smooth symplectic algebraic variety, 
and the projection $\pi$ is a symplectic resolution of
singularities \cite[Theorem 12]{Maf}.  The variety $\tilde{S}^\mu_\nu$ deformation retracts onto
the Spaltenstein variety $X_\nu^\mu$, and the irreducible components of
$X_\nu^\mu$ are Lagrangian in $\tilde{S}^\mu_\nu$.
We conjecture that $A^\mu_\nu$ is isomorphic to the Ext-algebra of a
sum of Lagrangians in the Fukaya category of $\tilde{S}^\mu_\nu$.  This conjecture is
completely analogous to the one that we made above involving $B(\cV)$
and the hypertoric variety $\MV$ in Remark \ref{hyp-conj}.  
When $N_\mu$ is the zero matrix, so that $\tilde{S}^\mu_\nu$ is
isomorphic to the cotangent bundle of the partial flag variety, our conjecture follows from the
Beilinson-Bernstein localization theorem \cite{BB} and the work of
Kapustin, Witten, Nadler, and Zaslow relating the Fukaya category of a cotangent bundle
to the category of perverse sheaves on the base \cite{KW07,NZ}.
\end{remark}

\section{Deformed category $\cO$}\label{sec:deformed O}
In this section we prove Theorem \ref{cat O and Spaltensteins}.
In order to understand the universal deformation $\tilde{A}_\nu^\mu$,
we will compare it to a ring $\wh{A}_\mu^\nu$ coming from the
``deformed category $\cO$'' considered by
Soergel \cite{Soe90,Soe92} and Fiebig 
\cite{Fie03, Fie06, Fie08}.  Here the deformation
comes from deforming the action of the Cartan
subalgebra.  Results of Fiebig and Soergel
allow us to show that $\wh{A}_\mu^\nu$ carries
a formal grading in the sense of Remark \ref{universal remark}, so it comes from the universal
deformation by extension of scalars.  It is
easy to construct deformed standard objects in 
the deformed category $\cO$, and to compute their
central characters.  We use this to show that 
the ``formal localization algebra" $\cZ(\wh{A}_\mu^\nu)$ is isomorphic to the completion of the 
the equivariant cohomology ring of the Spaltenstein $X^\mu_\nu$, but for the
larger torus $T^\mu\supset T^\mu_\nu$.

Most of the following material on deformed category
$\cO$ can be found in Fiebig's paper \cite{Fie03}.
However, he does not treat the parabolic case; 
when necessary, we will indicate how his arguments can be extended;
see also \cite[\S 2]{Str09}.
Let $\mu$, $\nu$ be compositions of $n$ as in the previous 
section, and assume that $\mu^+ \le \nu^t$, so 
the block $\cO^\mu_\nu$ is nonzero.  
Let $S^\mu := \Sym(\mt^\mu)^* \cong \Sym\mh^\mu$.  
A {\bf \boldmath$\mu$-deformation
algebra} is a commutative noetherian $S^\mu$-algebra $D$
with structure map
$\tau\colon S^\mu \to D$.  Given a 
$\mu$-deformation algebra and a weight $\alpha \in \Lambda\subset \mh^*$,
the {\bf \boldmath$\alpha$-weight space} of a left $U(\mg) \otimes D$-module
$M$ is
\[M_\alpha := \{v \in M \mid Hv = (\alpha(H) + \tau(H^\mu))v \;\text{for all}\; H \in \mh\},\]
where $H^\mu$ denotes the image of $H$ in the quotient $\mh^\mu$ of $\mh$.
Let $\mp = \mp_\mu$ be the parabolic subalgebra determined by $\mu$, 
as defined in the previous section, and let $D$ be a $\mu$-deformation algebra.  

\begin{definition}
The $D$-deformed 
$\mu$-parabolic category $\cO^\mu_{\! D}$ is the category
of finitely generated $U(\mg) \otimes D$-modules
$M$ such that
\begin{itemize}
\item $M = \bigoplus_{\alpha \in \Lambda} M_\alpha$,\,\,\,\,\, and
\item $(U(\mp) \otimes D)v$ is a finitely generated $D$-module for all $v \in M$.
\end{itemize}
For instance, if $D = \C$ and the map $S^\mu \to \C$ kills
$\mh^\mu$, then $\cO^\mu_\C$ is the usual parabolic
category $\cO^\mu$.
If $D \to D'$ is a map to another $\mu$-deformation algebra,
then $M \mapsto M \otimes_D D'$ defines a base change
functor $\cO^\mu_{\! D} \to \cO^\mu_{D'}$.  
\end{definition}

Deformed standard objects in $\cO^\mu_{\! D}$, which we call {\bf deformed Verma modules}, are defined in the following way.
Let $\mm = Z_\mg(\mt^\mu) \subset \mp$ be the centralizer 
of $\mt^\mu$; it is the subalgebra of 
block diagonal matrices for the $\mu$-blocks.  
For each $\alpha \in \Lplusmu$
there is a finite-dimensional irreducible $\mm$-module
$E_\alpha$ with highest weight $\alpha - \rho$, and these
give all finite-dimensional irreducible modules with integral weights.
For $\alpha \in \Lplusmu$, we define a deformed Verma module 
\[M^\mu_D(\alpha) := U(\mg) \otimes_{U(\mp)} (E_\alpha \otimes D),\]
where $U(\mp)$ acts on $E_\alpha$ via the map
$U(\mp) \to U(\mm)$ obtained by projecting away the off-diagonal
blocks of $\mp$, and on $D$ via 
$U(\mp) \to S^\mu \stackrel{\tau}{\to} D$.  It is
an object of $\cO^\mu_{\! D}$, where $D$ acts only on the last
factor.  Since it is generated as a module over $U(\mg) \otimes D$
by $1 \otimes v \otimes 1$, where $v \in E_\alpha$ is a
highest weight vector, the action of $D$ induces an isomorphism
$\End_{\wh{\cO}^\mu}(\wh{M}^\mu_D(\alpha)) \cong D$.
For any map $D \to D'$ of $\mu$-deformation
algebras, we have a natural isomorphism 
$M_D^\mu(\alpha)\otimes_D D' = M_{D'}^\mu(\alpha)$.
The object $M^\mu_\C(\alpha)$ is the standard cover of the simple module
$L(\alpha)$ in the usual parabolic category $\cO^\mu$.

To apply Fiebig's results we need our deformation algebra to be local.  
Let $\wh{S}^\mu = \prod_{i \ge 0} S^\mu_i$ be the 
completion of $S^\mu$ at the graded maximal ideal $S^\mu_{>0}$, and let
$\wh{\cO}^\mu$ denote the corresponding deformed category
$\cO^\mu_{\wh{S}^\mu}$.  We denote the deformed Verma modules in 
this category by $\wh{M}^\mu(\alpha) := M^\mu_{\wh{S}^\mu}(\alpha)$.

\begin{definition}
A {\bf Verma flag} for an object $M\in\wh{\cO}^\mu$ is a
finite filtration with subquotients that are isomorphic to deformed Verma modules.
\end{definition}

\begin{theorem}\label{deformed O}{\em\cite[\S 2]{Fie03}}
The category $\wh{\cO}^\mu$ has enough projectives.  The base change functor
$\wh{\cO}^\mu \to \cO^\mu_\C = \cO^\mu$ induces 
bijections between isomorphism classes of 
simples in $\wh{\cO}^\mu$ and in $\cO^\mu$, and
between isomorphism classes of indecomposable 
projectives in both categories.
All projective objects in $\wh{\cO}^\mu$ have 
Verma flags.
If $P \in \wh{\cO}^\mu$ is projective and $M\in\wh{\cO}^\mu$ has a Verma flag, then
$\Hom_{\wh{\cO}^\mu}(P, M)$ is a free $\wh{S}^\mu$-module, and
the natural map
\[\Hom_{\wh{\cO}^\mu}(P, M) \otimes_{\wh{S}^\mu} \C \to \Hom_{\cO^\mu}(P \otimes_{\wh{S}^\mu} \C, M\otimes_{\wh{S}^\mu} \C)\] 
is an isomorphism.
\end{theorem}

\begin{remark}
Note that Fiebig does not treat the parabolic case.  He assumes that
$\mp$ is a Borel subgroup, meaning that $\mu_i \le 1$ for all $i$,
so $W_\mu = \{1\}$.  We call such a composition {\bf regular}.
The main argument in \cite{Fie03} that needs modifying when 
$\mu$ is not regular is Lemma 2.3, which constructs
the projectives.  He explains that arguments of Rocha-Caridi and Wallach
\cite{RCW} can be adapted to the deformed situation.  The arguments
of \cite{RCW} do cover the parabolic case, so extending Fiebig's 
arguments is straightforward.
\end{remark}

Theorem \ref{deformed O} gives us for
each $\alpha \in \Lplusmu$
a projective object $\wh{P}^\mu(\alpha)$ such that 
$\wh{P}^\mu(\alpha)\otimes_{\wh{S}^\mu} \C \cong P^\mu(\alpha)$.  
If $\nu$ is another composition of $n$, we define $\alpha = \alpha_\nu$
as before, and let $\wh\cO^\mu_\nu$ be 
the full subcategory of $\wh\cO^\mu$ whose objects are all quotients
of direct sums of $\wh{P}^\mu(w\alpha)$ for $w \in \cI_\mu^\nu$.  Then
$\wh{P}^\mu_\nu := \bigoplus_{w \in \cI^\mu_\nu} \wh{P}^\mu(w\alpha)$ is 
a projective generator of this category, so $\Hom_{\wh{\cO}^\mu_\nu}(\wh{P}^\mu_\nu, -)$
defines an equivalence of categories
between $\wh\cO^\mu_\nu$ and finitely generated right modules over
$\wh{A}^\mu_\nu := \End(\wh{P}^\mu_\nu)$.  Theorem \ref{deformed O} also 
implies that the image of $\wh{M}^\mu(w\alpha_\nu)$ satisfies the formal analogues of the
hypotheses of Corollary \ref{deformed standard characterization}, so 
it is isomorphic to the deformed standard object in the category 
of $\wh{A}^\mu_\nu$-modules which we defined in Section \ref{sec-malleable}.

The base change functor
$\wh\cO^\mu_\nu \to \cO^\mu_\nu$ 
sends $\wh{P}^\mu_\nu$ to $P^\mu_\nu$, so it induces a ring homomorphism
$\wh{A}^\mu_\nu \to A^\mu_\nu$.  Theorem \ref{deformed O} implies
that $\wh{A}^\mu_\nu$ is a flat deformation of $A^\mu_\nu$ over $\Spec \wh{S}^\mu$.  
We wish to use Remark \ref{universal remark}
to relate this deformation to the universal deformation $\tilde{A}^\mu_\nu$.
To do this, we need to construct a formal grading on $\wh{A}^\mu_\nu$.

When $\mu$ is regular, we use a geometric interpretation
of deformed category $\cO$ due to Soergel \cite{Soe92} and 
Fiebig \cite{Fie03, Fie08} to construct our 
formal grading; we then deduce the case when $\mu$ is
general from this.  So assume for the moment that
$\mu$ is regular.  To indicate this, we omit the
superscript $\mu$ from our notations.  

As in the last section, let $\Fl_\nu$ denote the
partial flag variety $G/P_\nu$, and let $T\subset G$ be the diagonal subtorus acting on $X_\nu$.
Let $S := \Sym\mathfrak{t}^* \cong H_T^*(pt)$, and let $\wh{S} := \prod_{i=1}^\infty S_i$ be the completion
of $S$ at the graded maximal ideal.
For an element $w \in W/W_\nu$, let $C_w \subset \Fl_\nu$ denote the
$B$-orbit containing the $T$-fixed point $p_w$.
Let $\wh Z$ be the center of $\wh A_\nu$.

\begin{theorem}\label{Fiebig's results}
We have $\wh S$-algebra isomorphisms
$$H^*_T(\Fl_\nu) \otimes_S \wh{S} \,\,\,\cong\,\,\, \wh{Z} \,\,\,\cong\,\,\,
\End_{\wh{\cO}_\nu}(\wh{P}_0),$$
where $\wh{P}_0 := \wh{P}(w_0\alpha_\nu)$ is the antidominant projective.
The functor
\[\V = \Hom_{\wh{\cO}_\nu}(\wh{P}_0,-)\colon \wh{\cO}_\nu \to \wh{Z}-\mathrm{mod}\]
is full and faithful on objects with a Verma flag (in particular, on projectives) 
and we have natural $\wh{Z}$-module isomorphisms
$$\V\wh{P}(w\alpha_\nu) \cong \IH^*_T(\overline{C_w}) \otimes_S \wh{S}
\,\,\,\,\,\,\text{and}\,\,\,\,\,\,
\V\wh{M}(w\alpha_\nu) \cong H^*_T(C_w) \otimes_S \wh{S}$$
for all $w \in W/W_\nu$.
\end{theorem}

\begin{proof} The identification of the center of $\wh{A}_\nu$ is 
accomplished in \cite[Theorem 9]{Soe90} and \cite[3.6]{Fie03}.  Fiebig's
proof is instructive from our point of view: he shows that the map
\[\wh{Z} \to \bigoplus_{w\in W/W_\nu} \End_{\wh{\cO}_\nu}\!\left(\wh{M}(w\alpha)\right) 
\,\,\,\cong \bigoplus_{w\in W/W_\nu} \!\!\wh{S} \,\,\,\cong\,\,\, 
H^*_T(X_\nu^T) \otimes_S \wh{S}\]
is an injection, and the relations cutting out the image are the same as those 
that describe the image $H^*_T(X_\nu) \to H^*_T(X_\nu^T)$ in terms of the $T$-invariant curves in
$X_\nu$ \cite[1.2.2]{GKM}.  Since tensoring with $\wh{S}$ is exact for graded
$S$-modules, this gives the first isomorphism in Theorem \ref{Fiebig's results}.
It also gives the identification of $\V\wh{M}(w\alpha_\nu)$
with $H^*_T(C_w) \otimes_S \wh{S}$.

The full faithfulness of $\V$ is proven in \cite[7.1]{Fie08}.  Finally, 
\cite[7.6]{Fie08} identifies $\V\wh{P}(w\alpha_\nu)$
with the completion of sections of a sheaf on a ``moment graph" constructed
from the zero and one-dimensional orbits of $X_\nu$.  By
\cite[\S 2]{BrM01} this gives exactly $\IH^*_T(\overline{C_w}) \otimes_S \wh{S}$. 
\end{proof}

\begin{corollary} We have $\wh{S}$-algebra isomorphisms
\[\wh{A}_\nu \,\,\,\,\cong\,\,\,\,  \End_{\wh{Z}}\!\left(\bigoplus_{w \in W/W_\nu} \!\!\V\wh{P}(w\alpha_\nu)\right) \,\,\cong\,\,\,\, 
 \End_{H^*_T(\Fl)}\!\left(\bigoplus_{w \in W/W_\nu} \!\!\IH^*_T(\overline{C_w}) \right)\otimes_S \wh{S}.\]
\end{corollary}

This gives our formal grading of $\wh{A}_\nu$ in the non-parabolic case: the 
$i^\text{th}$ graded piece consists of maps which increase degree of the intersection
cohomology groups on the right by $i$.  
The grading in the general case now arises from the following proposition.

\begin{proposition}\label{deformed truncation}
Let $\mu, \nu$ be arbitrary compositions of $n$.  There is a surjective map from $\wh{A}_\nu$
to $\wh{A}^\mu_\nu$, with kernel generated by 
the idempotents $\{e_w\mid w \in \cI_\nu \setminus \cI^\mu_\nu\}$ and
the kernel of the projection $\mh \to \mh^\mu$. 
\end{proposition}

\begin{remark}
In fact, the kernel is generated by the idempotents alone,
but we do not need this, and will not prove it. 
\end{remark}

\begin{proof}
We construct a truncation functor $\hat\tau\colon \wh{\cO}_\nu \to \wh{\cO}^\mu_\nu$
which is the deformed analogue of the functor that takes the maximal
$\mp$-locally finite quotient of objects in category $\cO$.  We do this in 
two steps.  First, for an object $M \in \wh\cO_\nu$, define $M^\mu := M \otimes_{\wh S} \wh{S}^\mu$, 
the image of $M$ under the functor $\wh{\cO} \to \cO_{\wh{S}^\mu}$.  Next, define
\[Q\,\, := \bigoplus_{w \in \cI_\nu \setminus \cI^\mu_\nu} \wh{P}(w \alpha)^\mu,\]
and define $\hat\tau M$ to be the cokernel of the natural map  
$\Hom_{\cO_{\wh{S}^\mu}}(Q, M^\mu) \otimes_{\wh{S}^\mu} Q \to M^\mu$.
It is the largest quotient of $M^\mu$ which contains no subquotients isomorphic
to any simple object $L(w\a)^\mu$ with $w \in \cI_\nu \setminus \cI^\mu_\nu$.

It follows that this functor does indeed send $\wh{\cO}_\nu$ to $\wh{\cO}^\mu_\nu$.
It is not hard to see that $\hat\tau$ is right exact and left adjoint to 
the inclusion $\iota\colon \wh{\cO}^\mu_\nu \to \wh{\cO}_\nu$, so it 
sends projectives to projectives, and in fact sends a projective 
generator of $\wh{\cO}_\nu$ to a projective generator of $\wh{\cO}^\mu_\nu$.
This gives a natural homomorphism 
$\wh{A}_\nu \to \wh{A}^\mu_\nu$, which clearly contains the ideal described in the statement of the theorem.  
It is also surjective, since the adjunction map
$M \to \iota\hat\tau M$ is surjective for any $M$.
Our description of the kernel now follows from the characterization
of $\hat\tau M$ in the previous paragraph.
\end{proof}

Let $\tilde A_\nu^\mu$ be the universal deformation of $A_\nu^\mu$.
Let $S^\mu_\nu$ be the symmetric algebra of $Z((A^\mu_\nu)^!)^*_2$, so $\tilde A_\nu^\mu$
is an $S^\mu_\nu$-algebra.
Theorem \ref{flat}, Remark \ref{universal remark}, and Proposition \ref{deformed truncation} now imply the following result.

\begin{theorem}\label{relating the deformations}
There exists a linear map $\psi^\mu_\nu\colon Z((A^\mu_\nu)^!)^*_2 \to \mh^\mu$, inducing a graded ring homomorphism
$S^\mu_\nu \to \wh{S}^\mu$, such that we have
an $\wh{S}^\mu$-algebra isomorphism
\[\wh{A}^\mu_\nu \,\,\,\cong\,\,\, \tilde{A}^\mu_\nu \otimes_{S_\nu^\mu} \wh{S}^\mu.\]
\end{theorem}

By itself, Theorem \ref{relating the deformations} doesn't help us understand the 
universal deformation.  For instance, if it turned out that
$\psi^\mu_\nu = 0$, then we would have
$\wh{A}^\mu_\nu \cong A^\mu_\nu \otimes \wh{S}^\mu$. 
However, the following result implies that 
$\wh{A}^\mu_\nu$ carries all the information of the 
universal deformation.

\begin{proposition} \label{full deformation}
 The map $\psi^\mu_\nu$ of Theorem \ref{relating the deformations} is injective.
As a result, we have an $S_\nu^\mu$-algebra isomorphism
\[\tilde{A}^\mu_\nu \,\,\,\cong\,\,\, \bigoplus_{i \ge 0} \,\left(\wh{A}^\mu_\nu\right)_{\! i} \otimes_{S^\mu} S^\mu_\nu,\]
where the map $S^\mu \to S^\mu_\nu$ comes from any left inverse of $\psi^\mu_\nu$. 
\end{proposition}

We postpone the proof of Proposition \ref{full deformation} until we have established some further properties
of the deformation $\wh{A}^\mu_\nu$.  Note that since the maps 
$\wh{A}_\nu \to \wh{A}^\mu_\nu \to A^\mu_\nu$ are surjective, they induce
maps between the centers of these algebras.

\begin{lemma}\label{gln surjectivity}
The maps $Z(\wh{A}_\nu) \to Z(\wh{A}^\mu_\nu) \to Z(A^\mu_\nu)$ are surjective.
\end{lemma}
\begin{proof}
 By \cite[Theorem 2]{Bru08}, the action of the center of 
the enveloping algebra induces a surjection 
$Z(U(\mg)) \to Z(A^\mu_\nu)$.  It follows that 
$Z(U(\mg)) \otimes \wh{S}^\mu$ surjects onto $Z(\wh{A}^\mu_\nu)$ and (as a special case)
$Z(U(\mg))\otimes\wh{S}$ surjects onto $Z(\wh{A}_\nu)$.
The result then follows from the surjectivity of the maps 
$Z(U(\mg))\otimes\wh{S} \to Z(U(\mg)) \otimes \wh{S}^\mu \to Z(U(\mg))$.
\end{proof}

\begin{remark}
Proposition \ref{full deformation} and Lemma \ref{gln surjectivity} together imply that the universal deformation
$\tilde{A}^\mu_\nu$ is flexible in the sense
of Section \ref{sec:flexible-algebr}.
\end{remark}

The center of $\wh{A}^\mu_\nu$ acts on the deformed Verma module $\wh{M}^\mu(w\alpha)$ by a character
\[h_{w}^\mu \colon Z(\wh{A}^\mu_\nu) \to \End_{\cO^\nu_\nu}(\wh{M}^\mu(w\alpha)) \cong \wh{S}^\mu.\]
When $\mu$ is regular, which we indicate as usual 
by omitting $\mu$ from the notation, 
Theorem \ref{Fiebig's results} identifies $h_{w}$ with the map
$H^*_T(\Fl_\nu) \otimes_S \wh{S} \to H^*_T(p_w) \otimes_S \wh{S} \,\cong\, \wh{S}$ induced by restriction.
In terms of the identification
\[H^*_T(\Fl_\nu) \,\,\,\cong\,\,\, H_G^*(X_\nu)\otimes_{H_G^*(pt)}H_T^*(pt)
\,\,\,\cong\,\,\, H_{P_\nu}^*(pt)\otimes_{H_G^*(pt)}H_T^*(pt)
\,\,\,\cong\,\,\, S^{W_\nu} \otimes_{S^W} S,\]
we have $h_w(f \otimes g) = g \cdot w(f)$.
In particular, in degree 2 we have 
\begin{equation}\label{characters}
H^2_T(\Fl_\nu) \cong \mh^\nu \oplus \mh\,\,\,\,\,\text{and}\,\,\,\,\, h_w(x, y) = y + w(x).
\end{equation}
The characters of deformed parabolic Verma modules are determined by
the non-parabolic characters via the following result.

\begin{lemma}\label{parabolic characters}
 The diagram
\[\xymatrix@C=1.2cm{Z(\wh{A}_\nu) \ar[d]\ar[r]^(.57){h_{w}} & \wh{S} \ar[d] \\
Z(\wh{A}^\mu_\nu) \ar[r]^(.57){h^\mu_{w}} & \wh{S}^\mu 
}\]
commutes, where the left vertical map is the natural projection.
\end{lemma}

\begin{proof}
Using the proof of Proposition \ref{deformed truncation}, we see that the
deformed parabolic Verma module $\wh{M}^\mu(w\alpha) \cong \hat{\tau}\wh{M}(w\a)$ is 
a quotient of $\wh{M}(w\alpha)\otimes_{\wh{S}}\wh{S}^\mu$.
\end{proof}

\begin{remark}
Proposition \ref{full deformation}, Lemma \ref{parabolic characters}, and our formula for the characters $h_w$
in Equation \eqref{characters} together imply that the universal deformation
$\tilde{A}^\mu_\nu$ is malleable in the sense of Section \ref{sec-malleable}.
\end{remark}

We now use these calculations to determine the center of $\wh{A}^\mu_\nu$.
This result, along with Proposition \ref{full deformation}, will complete the proof of Theorem \ref{cat O and Spaltensteins}.

\begin{theorem}  The formal localization algebra
$\cZ(\wh{A}^\mu_\nu)$ is isomorphic to $H^*_{T^\mu}(X^\mu_\nu)\otimes \wh{S}^\mu$.
\end{theorem}

\begin{proof}
Consider the commutative diagram
\[\xymatrix@R=0.8cm@C=1cm{Z(\wh{A}_\nu) \ar[d]\ar[rr]^(.45){\oplus h_{w}} & & \bigoplus\limits_{w\in \cI_\nu}\wh{S} \ar[d]^q \\
Z(\wh{A}^\mu_\nu) \ar[rr]^(.45){\oplus h^\mu_{w}} & & \bigoplus\limits_{w\in\cI^\mu_\nu}\wh{S}^\mu 
}\]
where $q$ is the quotient map $\wh{S} \to \wh{S}^\mu$ for all $w \in \cI^\mu_\nu$, and kills all
terms for $w \notin \cI^\mu_\nu$.  The horizontal maps are injective by 
Equation \eqref{characters} and Lemma \ref{parabolic characters}, while the left vertical map is surjective by Lemma \ref{gln surjectivity}. 
Thus we have an isomorphism $Z(\wh{A}^\mu_\nu) \cong q(\Im \oplus h_{w})$.

This diagram has a topological analogue:
\[\xymatrix@R=0.8cm@C=1cm{H^*_T(\Fl_\nu) \otimes_S \wh{S} \ar[d]\ar[rr]^(.55){\iota^*} & & \bigoplus\limits_{w\in \cI_\nu}\wh{S} \ar[d]^q \\
H^*_{T^\mu}(\Fl^\mu_\nu) \otimes_{S^\mu} \wh{S}^\mu \ar[rr]^(.6){(\iota^\mu)^*} & & \bigoplus\limits_{w\in\cI^\mu_\nu}\wh{S}^\mu 
}\]
 The left vertical map is the composition of the 
restriction to the subtorus $T^\mu$ with the restriction to $\Fl^\mu$.
The maps $\iota^*$ and $(\iota^\mu)^*$ are the restrictions to the fixed point sets, which are indexed as indicated by Proposition 
\ref{Spaltenstein fixed points}.

Both $\Fl_\nu$ and $\Fl_\nu^\mu$ have vanishing odd cohomology, so they are equivariantly formal, which 
implies that the horizontal maps are injections.  The restriction from
$H^*(\Fl_\nu)$ to $H^*(\Fl_\nu^\mu)$ is surjective by \cite[Corollary 2.5]{BrO}.
Together with equivariant formality this implies that the 
left vertical map is surjective, so we have an isomorphism 
$H^*_{T^\mu}(\Fl^\mu_\nu) \cong q(\Im \iota^*)$.  We have already seen in the proof of Proposition \ref{Fiebig's results}
that $\iota^*$ is identified
with $\oplus h_w$, so this proves the theorem.
\end{proof}

By Lemma \ref{parabolic characters} the image of $(x, y)$ under the composition
\[\mh^\nu \oplus \mh \cong Z(\wh{A}_\nu)_2 \to Z(\wh{A}^\mu_\nu)_2 \to \End(\wh{M}^\mu(w\alpha))\]
is multiplication by $(y + w(x))^\mu \in \mh^\mu \subset \wh{S}^\mu$.
It follows from \cite[Theorem 2]{Bru08} that any two
simples in $\cO^\mu_\nu$ can be connected by a chain of
non-trivial $\operatorname{Ext}^1$ groups, so by Lemma \ref{image of psi} we see that 
the image of $\psi^\mu_\nu$ contains 
\[\Sigma := \{(w(\alpha_\nu) - v(\alpha_\nu))^\mu \mid v, w \in \cI^\mu_\nu\}\subset \mh^\mu.\]
Using this, we can finally prove Proposition \ref{full deformation}.

\begin{proof}[Proof of Proposition \ref{full deformation}]
We have  $\dim \Span(\Sigma) \le \rank \psi_\nu^\mu \le \dim Z(A^{\nu^o}_{\mu})_2$, 
so it will be enough to show that $\dim \Span(\Sigma) = \dim Z(A^{\nu^o}_{\mu})_2$.

As we noted in Remark \ref{Brundan}, 
Brundan \cite{Bru06} shows that $Z(A^{\nu^o}_{\mu})_2 \cong \mt^\mu_\nu$.
Lemma \ref{extra restrictions} shows that 
the pairing between $\mh^\mu$ and $\mt^\mu$ induces 
a well-defined pairing between $\Sigma$  (a subspace of $\mh^\mu$) and
$\mt^\mu_\nu$ (a quotient of $\mt^\mu$).  We will show that this 
pairing is non-degenerate.

As we saw in the proof of Lemma \ref{extra restrictions}, the 
set $\{w(\alpha_\nu) \mid w \in \cI^\mu_\nu\}$  is the set of 
all vectors in $\Z^n$ which have 
$\nu_{-i}$ entries equal to $i$ for all $i \in \Z$, and for which  
the entries in each $\mu$-block are strictly decreasing.
It follows that if the $\mu$-blocks are reordered, the 
effect on the set $\Sigma$ is just to apply the appropriate
permutation to each element.  The same holds for $\mt^\mu_\nu$, 
so without loss of generality we can assume that $\mu = \mu^+$.  It is
also easy to see that we can take $\nu = \nu^+$.

We must show that every element in $\mt^\mu$ which 
pairs to zero with $\Sigma$ must lie in the span of $\{\mathbf{1}_J\mid J \in \mathcal{J}\}$.  
To do this, let $k_1 < \dots < k_r$ be
the solutions $k$ to the equation
\[\mu_1 + \dots + \mu_k = \lambda_1 + \dots + \lambda_k,\]
and let $J_i$ be the union of the first $k_i$ 
$\mu$-blocks.  Then Lemma \ref{extra restrictions}
says that for any vector $w \alpha_\nu$ and any $i$ the multiset of entries
indexed by $j \in J_i \setminus J_{i-1}$ is independent of $w$. 
It is also clear that the entries in each $J_i \setminus J_{i-1}$ 
can be chosen independently.  It follows that without loss of
generality we can assume that $\mathcal{J} = \left\{ \{1, \dots, n\} \right\}$.

We make this assumption, and proceed by induction on the number of
nonzero entries in $\mu$. If there is only one, then $\mt^\mu_\nu = 0$
and we are done.  Otherwise, consider filling the $\mu$-blocks with
entries of $\alpha_\nu$, starting with the left-most block first.  If
we fill the first block with the entries $1,\dots,\mu_1$, then the
remaining blocks give an element of $\{w(\alpha_{\nu'}) \mid w \in
\cI^{\mu'}_{\nu'}\}$ for the pair
\[\mu' = (\mu_2, \mu_3, \dots) \;\text{ and }\; \nu' = (\nu_1 - 1, \dots,
\nu_{\mu_1} - 1, \nu_{\mu_1 + 1}, \dots)^+.\]
We can describe the transpose partition $\lambda' := (\nu')^t$ 
as follows: if $m$ is the unique integer such that $\lambda_m \ge \mu_1$ and $\lambda_{m+1} < \mu_1$,
then 
\[\lambda' = (\lambda_1, \dots, \lambda_{m-1}, \lambda_{m+1} + \lambda_{m} - \mu_1, \lambda_{m+2}, \dots).\] 
Note that for these new partitions, we have $\mu' < \lambda'$ and
$\mathcal{J}=\{\{1,\dots,n\}\}$
\begin{equation*}
   (\lambda_1' + \dots + \lambda_k')=  (\lambda_1 + \dots + \lambda_k) >
    (\mu_1 + \dots + \mu_k)  \geq  (\mu_2 + \dots + \mu_{k+1}) =  (\mu_1' + \dots + \mu_k')
    \,\,\,\,\,\,\text{if $k<m$}
\end{equation*}
and
\begin{equation*}
   (\lambda_1' + \dots + \lambda_k') = (\lambda_1 + \dots + \lambda_{k+1}-\mu_1) > (\mu_2 \dots + \mu_{k+1}) = (\mu_1' + \dots + \mu_k') 
   \,\,\,\,\,\,\text{if $k\geq m$.}
\end{equation*}

Using permutations that keep the first block fixed, we
obtain an inclusion $\Sigma^{\mu'}_{\nu'}\hookrightarrow
\Sigma^{\mu}_{\nu}$, and by our inductive hypothesis, no element of
the span of this subset is annihilated by $\mt^{\mu}_{\nu}$, which
pairs via the surjective quotient map $\mt^{\mu}_{\nu}\to \mt^{\mu'}_{\nu'}$.

The kernel of this map is one dimensional, spanned by the element $\mathbf{1}_{[1,\mu_1]}$
which has all 1's on the first block and 0's elsewhere.  
So by the inductive hypothesis, we only need to find an element of
$\Sigma$ which pairs non-trivially with this vector. 
That is, we must find $v$ and $w$ in $\cI^\mu_\nu$ for which 
 $v(\alpha_\nu)$ and $w(\alpha_\nu)$ have different entries 
in the first block.  Rather
than construct the whole vector, we note that a choice of entries in
the first block can be extended to a vector of the form $w(\alpha_\nu)$
if and only if the remaining partitions $\mu''$ and $\nu''$ still
satisfy the dominance condition $(\mu'')^+\leq \la''$.

We have already noted that we can take the entries in our first
block to be $1,\dots,\mu_1$.  We claim that $1,\dots,\mu_1-1,\mu_1+1$
will also extend to a vector of the form $v(\alpha_\nu)$ for $v\in \cI^\mu_\nu$.  This will
finish the proof, since the difference of these vectors will pair
non-trivially with $\mathbf{1}_{[1,\mu_1]}$, and thus establish
non-degeneracy.

In this case, we have 
\[\nu'' = (\nu_1 - 1, \dots, \nu_{\mu_1 - 1} -1, \nu_{\mu_1}, \nu_{\mu_1 + 1} - 1, \nu_{\mu_1 + 2}, \dots)^+.\]
If $\nu_{\mu_1} = \nu_{\mu_1 + 1}$, then $\nu' = \nu''$, so
$(\mu'')^+\leq \la''=\la'$ and we are done.  
Otherwise, we find that
\[(\nu'')^t = (\lambda_1, \dots, \lambda_{m-2}, \lambda_{m-1} - 1, \lambda_{m+1} + \lambda_{m} - \mu_1 + 1, \lambda_{m+2}, \dots).\]
Since $\lambda_1 + \dots + \lambda_{m-1} > \mu_1 + \dots + \mu_{m-1}$ by assumption, we have
$\mu \le (\nu'')^t$, and so it is possible to continue filling the remaining blocks.
\end{proof}
\end{spacing}

 \bibliography{./symplectic}
\bibliographystyle{amsalpha}
\end{document}